\documentclass[leqno,12pt]{article}
\usepackage[hang]{subfigure}
\usepackage{color}
\usepackage{multirow}
\usepackage{latexsym}
\usepackage{float}
\usepackage{hyperref}

\makeatletter
\newsavebox\myboxA
\newsavebox\myboxB
\newlength\mylenA

\newcommand*\xoverline[2][0.75]{%
    \sbox{\myboxA}{$\m@th#2$}%
    \setbox\myboxB\null
    \ht\myboxB=\ht\myboxA%
    \dp\myboxB=\dp\myboxA%
    \wd\myboxB=#1\wd\myboxA
    \sbox\myboxB{$\m@th\overline{\copy\myboxB}$}
    \setlength\mylenA{\the\wd\myboxA}
    \addtolength\mylenA{-\the\wd\myboxB}%
    \ifdim\wd\myboxB<\wd\myboxA%
       \rlap{\hskip 0.5\mylenA\usebox\myboxB}{\usebox\myboxA}%
    \else
        \hskip -0.5\mylenA\rlap{\usebox\myboxA}{\hskip 0.5\mylenA\usebox\myboxB}%
    \fi}
\makeatother

\newcommand{\var}{\mbox{Var}}
\newcommand{\E}{\mathbb{E}}
\newcommand{\cov}{\mbox{Cov}}

\newcommand{\R}{\mathbb R}

\newcommand{\N}{\mathbb{N}}

\def\P{\mathbb{P}}

\usepackage{verbatim}
\usepackage{amsmath,stackrel}
\usepackage{amssymb}
\usepackage{ntheorem}
\usepackage[ansinew]{inputenc}
\usepackage{graphicx}
\usepackage{epsfig}
\usepackage{dsfont}
\usepackage{bm}
\usepackage{bbm}

\usepackage[authoryear]{natbib}
\newtheorem{theorem}{Theorem}[section]
\newtheorem{prop}{Proposition}[section]

\newtheorem{assumption}{Assumption}[section]

\begin{document}

\title{Detecting relevant differences in the covariance operators of functional time series  - a sup-norm approach}

\author{
{\small Holger Dette, Kevin Kokot} \\
{\small Ruhr-Universit\"at Bochum} \\
{\small Fakult\"at f\"ur Mathematik}\\
{\small Bochum, Germany} \\
{\small e-mail: $\{$holger.dette, kevin.kokot$\}$@rub.de}\\
}

  \maketitle

\begin{abstract}
In this paper  we propose   statistical inference tools for the covariance operators 
of functional time series in the two sample and change point problem.
In contrast to most of the literature  the focus of our approach 
is not testing the null hypothesis of exact equality of the covariance operators. Instead we propose to formulate the null hypotheses in the form that ``the distance between the operators is small'', where we measure deviations  by the sup-norm.   We provide powerful  bootstrap  tests for these type of hypotheses, investigate their asymptotic properties  and study their finite sample properties by means of a simulation study.

 \end{abstract}

Keywords: covariance operator, functional time series, two sample problems, change point problems, CUSUM, relevant hypotheses, Banach spaces, bootstrap
\\ 
AMS Subject Classification:  62G10,  62M10

\section{Introduction}  \label{sec:Intro}
 \def\theequation{1.\arabic{equation}}
\setcounter{equation}{0}

The field of functional data analysis has
found considerable  attention in the  statistical literature
 as in many applications the observed data points exhibit certain degrees of
dependence and smoothness and thus may naturally be regarded as discretized
functions. Introductions to this topic  can be found in the monographs of
\cite{bosq2000}, \cite{RamsaySilverman2005},
\cite{FerratyVieu2010}, \cite{HorvathKokoskza2012} and \cite{hsingeubank2015}, among others.
Interest may, for example, be in comparing characteristic
parameters of the random functions from  two different samples
(\textit{two sample problem}) or in investigating whether a certain parameter
of a   functional time series  remains stable over time
(\textit{change point problem}). In most cases the  considered parameters (such as the mean)
 are functions themselves, which makes the analysis of this type of problems
 challenging.  In the present paper we investigate the
second-order properties of a stationary functional time series which are contained in its
covariance operators and  important for the
understanding of the smoothness of the  stochastic fluctuations of the data
\citep{Kraus2012}. Most
of the literature on this topic considers Hilbert space-valued random variables.
 The popularity  of this  approach is due
 to the fact  that  such a framework allows the development of
dimension reduction techniques such as (functional) principal components. On the other hand
dimension reduction may  yield to a loss of information
 as data is projected on  finite dimensional spaces, and several
 authors  argue that it might be more reasonable to work in the space of functions directly
 \citep[see, for example,][for a recent reference]{Aue2015DatingSB}.

In this paper we will develop methodology   to  compare the covariance operators of two functional time series and
 to detect changes in the covariance operator of a functional time series in the space of
continuous functions defined on a compact interval. Thus - in contrast to most
of the literature on this topic, which considers Hilbert space-valued objects -
 the random variables under consideration
 are (dependent) elements    of  a Banach space, and it is possible to compare the covariance operators in the sup-norm.
  Another important
 difference to the literature    consists in the fact
 that the main focus of our approach is not on {\it classical} hypotheses of the form
\begin{align} \label{classhyp}
H_{0}: C_1 = C_2  \qquad \mbox{versus} \qquad \ H_{1}: C_1 \not = C_2
\end{align}
 where $C_1$ and $  C_2$ are  either the   covariance operators corresponding
 to the  two samples or to the covariance operator  before and
after a change point. In contrast  we  consider {\it relevant} hypotheses of the form
\begin{align} \label{eq:relevant-intro}
  H_0^{ \Delta}: d(C_1, C_2) \leq \Delta \qquad \mbox{versus}
  \qquad H_1^{ \Delta}: d(C_1, C_2) > \Delta
\end{align}
where  $\Delta \geq  0 $ is a given  threshold and  $d$ a suitable metric on the space of covariance operators (in our case the $\sup$-norm).
Note that hypotheses of the form \eqref{eq:relevant-intro} contain the  classical hypotheses  in~\eqref{classhyp} as a special case for the choice
$\Delta =0$, but we argue that the  case $\Delta >0 $ is at least of equal interest. In fact, in many applications it is obvious  that   $C_1$ and $  C_2$  can not exactly coincide but the deviation
might be small.  In such cases testing for exact equality may be questionable and
it might be more reasonable to test for a relevant or  significant deviation between
the two covariance  operators.

In the case of testing classical hypotheses the metric does not matter  because
under the null hypothesis the  distance between $C_{1}$ and $C_{2}$
vanishes in any metric. However, this is not the case for relevant hypotheses
of the form \eqref{eq:relevant-intro}.
In the present  context  two covariance operators  with rather different shapes may still have a
small $L^2$-distance,  which makes an appropriate interpretation of the
threshold $\Delta$ for practitioners difficult. As an alternative  we propose to consider the maximum
deviation between  the covariance operators  as metric in the hypotheses
 \eqref{eq:relevant-intro}. On the one hand this metric
 makes the interpretation of  the threshold $\Delta$ more easy. On the other
hand it  leads to a Banach space-based framework where no dimension reduction
techniques are available and  the development and
theoretical justification of statistical
methods are more challenging.

In Section~\ref{sec:banach_methods} we review some basic properties of
random variables
in the space of continuous functions. In particular we  define moments of order two
through injective tensor products.
We also state a  central limit theorem for a stationary Banach-space valued process, which will be the basis for all theoretical arguments given in this paper.
 In Section~\ref{sec:TSP} we develop statistical methods
for the comparison of covariance operators in the two sample problem. In particular
a  test is proposed  for the null
hypothesis of no relevant difference between the covariance operators
from two independent samples.  As a special (and substantially simpler case)
we  also construct a new test for the classical hypotheses \eqref{classhyp} with a simple structure
and  nice statistical properties.
Section~\ref{sec4} is devoted to  the change point problem, where
methodology  is developed to detect changes in the covariance operator of a functional
time series. In all cases we make use of a multiplier  bootstrap procedure to obtain critical values
for the proposed tests. The theoretical justification
of all methods is given in Section~\ref{sec7}, while
Section~\ref{sec:simulations} contains a detailed simulation study
to investigate the finite sample properties of  the proposed
tests. Although classical hypotheses are not the main focus of our work  we also
 compare the new tests  for the classical hypotheses
with some of  the currently available methodology  and demonstrate that  they provide powerful
alternatives to the procedures, which have been proposed in the literature so far.

\subsection{Related literature}

There exists a considerable amount of  literature considering the comparison of covariance operators in the two sample
problem, where  random functions in the Hilbert space of  square-integrable functions
and the classical null hypothesis of equal covariance operators are  investigated.
\cite{Panaretos2010} consider independent
Gaussian data and describe an application to DNA minicircle data.
\cite{fremdt2013}  extend  the  theoretical findings of these authors to a  more general
 model such that non-Gaussian curves  are also covered. In both
references, functional principal components (FPCs) are used for dimension
reduction. \cite{Kraus2012} introduce the notion of a dispersion operator and
propose a robust test, which  is based on a  truncated
version of the Hilbert-Schmidt norm of a score operator defined via the
dispersion operator.  \cite{Zhang2015}  propose a pivotal test procedure
based on FPCs and  self-normalization and also provide inference tools for
the eigensystem of the covariance operators.

Several authors argue  that  dimension reduction may yield to a loss of information
and  propose alternative procedures  for the  comparison of covariance operators in the two sample
problem. \cite{pigoli2014} discuss different distance measures
between covariance operators and develop a permutation test and
\cite{Paparoditis2016}  propose a  bootstrap  test for the
(classical) null hypothesis of equality of $K$ covariance operators.
\cite{cabassi2017} suggest to combine all  pairwise comparisons between  samples of independent data into a global test for this problem, where the  Hilbert-Schmidt norm
between the square roots of the covariance operators is used as a measure of deviation.
 \cite{boente2018} provide a theoretical framework which clarifies the ability of the test
  to detect local alternatives.
 \cite{Pilavakis2019}   develop
 a fully functional  test for the equality of auto-covariance operators of
  temporally dependent time series, which is based on a moving
 block bootstrap.
For  independent data the $K$-sample problem has also been considered by
 \cite{Guo2016}
 who  propose to estimate
 the  supremum value of the sum of the squared differences between the estimated individual covariance functions
 and the pooled sample covariance function.

So far, the change point problem for covariance operators has found less
attention in the literature. \cite{jaruskova2013} uses FPCs to develop a test  for the existence of a change point, while
 \cite{stoehr2019}  use  the  circular block  bootstrap to construct a change point test.
In particular these authors  develop a test based on dimension reduction and two procedures
 which take the full functional structure into account.
 A fully functional test has also been proposed by
 \cite{sharipov2019}, who use a non overlapping block bootstrap to obtain critical values.
More recently,  \cite{Aue2020}  propose statistical tests
 for detecting a change in the spectrum
and in the trace of the covariance operator, respectively.

All these references consider the problem of testing classical hypotheses of the form
  \eqref{classhyp}. Recently  \cite{detkokvol2020} propose a   comparison
  of   covariance operators in the two sample problem and in the context of change point analysis  by
  testing relevant hypotheses of the form \eqref{eq:relevant-intro}, where an $L^2$-distance is used as metric.
 However, in the context of testing relevant hypotheses the norm  matters
 as two covariance operators might be close in one norm but not in another.
 In particular,   relevant deviations between covariance operators  in  the sup-norm have
 - to our best knowledge - not been considered so far
 and requires a different methodology as the space under consideration is a Banach but not
 a Hilbert space. There does not exist so much literature on  functional data analysis considering Banach spaces
 and exemplarily we mention the recent  work of  \cite{dette2018} who considered relevant hypotheses
 for the mean function  and \cite{lieblreim2019}, who developed confidence bands for functional parameters.

\section{$C(T)$-valued random variables} \label{sec:banach_methods}
\def\theequation{2.\arabic{equation}}
\setcounter{equation}{0}

In this paper  we consider random variables taking values in the Banach space
of real-valued and continuous functions defined on  a compact set $T$ and denote this space by $C(T)$,
which is equipped with the sup-norm $\|X \|_\infty =\max_{t \in T} |X(t)|$ for any $X \in C(T)$. The
underlying probability space $(\Omega, \mathcal{A}, \P)$ is assumed to be
complete and measurability is always meant with respect to the  Borel
$\sigma$-field  generated by the open sets relative to the respective
sup-norm.

Following Chapter~11 in \cite{jankai2015}, we use injective
tensor products to define moments of $C(T)$-valued random variables and note that
$C(T)^{\check \otimes k} = C(T^k)$  isometrically with the natural
identification (Theorem~11.6). The $k$th moment of  a $C(T)$-valued random variable
$X$  exists,
whenever $\E\big[\|X\|^k_\infty] < \infty$
(Theorem~11.25) and is defined by the function
in $C(T)^{\check \otimes k} = C(T^k)$, which maps $(t_1,\dots,t_k) \in T^{k}$ to
\begin{align*}
\E X^{\check \otimes k}(t_1,\dots,t_k) = \E \big[ X(t_1)\cdots X(t_k) \big]
\end{align*}
(Theorem~11.10). Throughout this paper, we write
$X^{\check \otimes 2} = X \check \otimes X$ for any $X\in C(T)$ and mean the function
in $C(T^2)$ defined by $(s,t) \mapsto X(s)X(t)$. Consequently, the covariance
operator of a $C(T)$-valued random variable is defined by
\begin{align*}
C(\cdot,\cdot) = \cov(X(\cdot),X(\cdot))
= \E \big[(X - \mu)^{\check \otimes 2}(\cdot,\cdot) \big] \in C(T^2)
\end{align*}
where $\mu = \E [X] \in C(T)$ is the expectation of $X$.

Let $\rho$ denote  a metric on $T$ such that $(T,\rho)$ is totally bounded, then
the metric  $\rho_{\max}$ on $T^2$ is defined through
$\rho_{\max} ((s,t),(s^\prime, t^\prime))  = \max\{\rho(s,s^\prime),\rho(t, t^\prime)\}$
and  the expression $D(\omega,\rho_{\max})$ denotes  the packing
number with respect to the metric $\rho_{\max}$ on $T^2$ that is the maximal
number of $\omega$-seperated points in $T^2$ \citep{wellner1996}.
Note that in this case
$(T^2,\rho_{\max})$ is totally bounded as well.

In order to describe the  dependence in the data
we introduce the   concept $\varphi$-mixing and
denote by $\mathbb{P}(G|F)$ the conditional probability of $G$ given $F$.
For  two $\sigma$-fields $\mathcal{F}$ and $\mathcal{G}$  we define the coefficient
\begin{align} \label{eq:mixing-coef}
\phi(\mathcal{F},\mathcal{G})
= \sup \big\{ |\mathbb{P}(G|F) - \mathbb{P}(G)| \colon
F\in \mathcal{F}, ~G\in\mathcal{G}, ~\mathbb{P}(F)>0 \big\} \, .
\end{align}
For a given strictly stationary sequence $(\eta_{j})_{j\in\mathbb{N}}$ of
random variables in $C(T)$, denote by $\mathcal{F}^{k^\prime}_{k}$ the
$\sigma$-field generated by $( \eta_{j}\colon k\leq j \leq k^\prime)$.
Then, the $k$th $\varphi$\textit{-mixing coefficient} of
$(\eta_{j})_{j\in\mathbb{N}}$ is defined by
\[
\varphi (k)
= \sup_{k^\prime \in \mathbb{N}} \phi (\mathcal{F}_{1}^{k^\prime} ,\mathcal{F}_{k^\prime+k}^\infty)
\]
and the stationary  time series $(\eta_{j})_{j\in\mathbb{N}}$ is
called $\varphi$\textit{-mixing} whenever the sequence of mixing coefficients
converges to zero as $k\to\infty$.

Given the preceding discussion, the analysis of the covariance operators of
random variables in $C(T)$ can in some sense be regarded to the analysis
 of $C(T^2)$-valued random variables. 
More  precisely, Theorem~11.7 in \cite{jankai2015} implies that $C(T^2)$ is separable
such  that
measurability issues are avoided, and Theorem~1.3 in \cite{billingsley1968}
implies that any $C(T^2)$-valued random variable is tight. A random function
$X$ in $C(T^2)$ is called Gaussian if and only if its finite dimensional
vectors $(X(t_1),\dots,X(t_k))$   follow a multivariate normal distribution for any $t_1,\dots,t_k \in T^2$ and $k\in\N$.

\begin{assumption} \label{assumption}
$(Z_j)_{j\in\mathbb{N}}$ is a sequence of $C(T)$-valued
random variables such that
\begin{align*}
  Z_{j} =  \mu + \eta_{j} \, , \quad j\in \N
\end{align*}
where $\mu \in C(T)$ denotes the expectation function and $(\eta_j)_{j\in\mathbb{N}}$ is a strictly stationary process.
\begin{enumerate}\itemsep-.4ex
\item[(A1)]
The packing number $D(\omega, \rho_{\max})$ satisfies
\begin{align*}
  \int_0^\tau  D(\omega,\rho_{\max})^{1/J}~d\omega <\infty
\end{align*}
for some $\tau > 0$ and some even integer $J \geq 2$.

\item[(A2)]
There is a constant $K$ such that
\begin{align*}
\mathbb{E}[\|\eta_{1}\|_\infty^{4+\nu}] \leq K \, ,
\quad \mathbb{E}[\|\eta_{1}\|_\infty^{2J}] \leq \infty
\end{align*}
for some $\nu>0$, where $J$ is the same integer as in (A1).
\item[(A3)]
There  exists a real-valued  non-negative
random variable $M$ and a constant $\tilde K$ such that, for any
$j\in\mathbb{N}$,
$\E \big[ (\|\eta_{j}\|_\infty \, M )^J \big] < \tilde K < \infty$ and the inequality
$$
|\eta_{j}(t)-\eta_{j}(t^\prime)|\leq  M \rho(t, t^\prime)
$$
holds almost surely for all $t,t^\prime\in T$. The integer $J$ is the same as in (A1).
\item[(A4)]
The process $(\eta_j)_{j\in\mathbb{N}}$ is $\varphi$-mixing with  mixing coefficients satisfying, for some $\bar\tau \in (   1/(2+2\nu) , 1/2 )$,  the condition
\begin{align*}
\sum_{k=1}^\infty k^{1/(1/2-\bar\tau) } \varphi(k)^{1/2}  < \infty \, , \quad
\sum_{k=1}^\infty (k+1)^{J/2-1}\varphi(k)^{1/J}  < \infty \, ,
\end{align*}
where the constants  $J$ and $\nu$ are the same  as in (A1) and (A2),
respectively.
\end{enumerate}
\end{assumption}

Note that
Assumption~\ref{assumption} implies  the existence of the covariance operator defined by
\begin{align} \label{h1}
C(s,t) = \cov(Z_{j}(s), Z_{j}(t)) = \E[(Z_{j}(s)-\mu(s))(Z_{j}(t)-\mu(t))] \, .
\end{align}
Condition (A4) on the summability of the mixing
coefficients is satisfied if there exists an $a \in (0,1)$ such that
$\varphi (k) \leq c a^k$  ($k\in \mathbb{N}$).
 For the formulation and a proof of a CLT of  Banach space valued random variables
 we denote by   the symbol ``$\rightsquigarrow$'' weak
convergence in $(C(T))^k $ or $(C(T^2))^k$ and the symbol
``$\stackrel{\mathcal{D}}{\longrightarrow}$'' denotes weak convergence in
$\R^k$ for some $k\in\N$.
The following  result  is proved in Section~\ref{sec7}.

\begin{theorem} \label{thrm:CLT}
Let $(Z_j)_{j\in\mathbb{N}}$ denote a stochastic process in $C(T)$ satisfying Assumption~\ref{assumption}.
Then,
\begin{align*}
G_n = \frac{1}{\sqrt{n}}
\sum_{j=1}^n ((Z_{j} - \bar{Z}_{n})^{\check \otimes 2} - C )
\rightsquigarrow Z
\end{align*}
in $C(T^2)$ as $n\to\infty$ where $\bar{Z}_{n} = 1/n\sum_{j=1}^n Z_{j} \in C(T)$, $C$ is defined by \eqref{h1} and $Z$ is a
centred Gaussian random variable with covariance operator 
\begin{align}\label{eq:long-run-cov}
\mathbb{C}((s,t),(s^\prime, t^\prime)) = \mathrm{Cov}(Z(s,t),Z(s^\prime, t^\prime))
= \sum_{i=-\infty}^{\infty} \cov\big(\eta_{0}^{\check \otimes 2}(s,t),
\eta_{i}^{\check \otimes 2}(s^\prime,t^\prime)\big) \, .
\end{align}
\end{theorem}

In the remaining part of the paper, we
consider the unit interval $T = [0,1]$ and, for a positive constant
$\theta \in (0,1]$, the metric $\rho(s,t) = |s-t|^\theta$ on $[0,1]$.
Consequently, on $T^2 = [0,1]^2$, we use the metric
$\rho_{\max} ((s,t),(s^\prime, t^\prime)) = \max\{\rho(s,s^\prime),\rho(t, t^\prime)\}$
and the packing number of the square $[0,1]^2$ with respect to this metric satisfies
$D(\omega,\rho_{\max})  \lesssim  \big\lceil \omega^{-2/\theta} \big\rceil $
(to see this, consider the points $(k\omega^{1/\theta}, l\omega^{1/\theta})$
for $k,l = 0,\dots,\lfloor\omega^{-1/\theta} \rfloor$). Therefore
condition (A1) reduces to
\begin{align*}
\int_0^\tau  D(\omega,\rho_{\max})^{1/J}~d\omega
&   \lesssim  \int_0^\tau \big\lceil \omega^{-2/\theta} \big\rceil^{1/J}~d\omega
\lesssim \frac{\tau^{1-2/(J\theta)}}{1-2/(J\theta)} < \infty
\end{align*}
and holds,
whenever the even integer $J$ satisfies $J>2/ \theta$ and consequently, under
this assumption, H\"older continuous processes satisfy (A1). Because the  paths of the
Brownian Motion $\{W(t)\}_{t\in[0,1]}$ are H\"older continuous of order
$\theta$ for any $\theta \in (0,1/2)$ and the random variable $\|W\|_\infty$ has moments of all order
 Assumption~\ref{assumption} is satisfied  for the Brownian motion (we can use $J=6$ in (A4) for  this case).
 For general processes with less smoothness, that is a smaller
constant $\theta$, we require a stronger summability assumption (A4) on
the mixing coefficients and the existence of higher moments.

\section{The two sample problem} \label{sec:TSP}
\def\theequation{3.\arabic{equation}}
\setcounter{equation}{0}

Throughout this section, we consider two independent samples
$(X_{j}\colon j=1,\dots,m)$ and $(Y_{j}\colon j=1,\dots,n)$ drawn from
independent strictly stationary sequences $(X_{j})_{j\in \N}$ and $(Y_{j})_{j\in \N}$ in $C([0,1])$ with representations
\begin{align} \label{mod}
  X_j =  \mu_1 + \eta_{1,j} \, , \quad Y_j = \mu_2 + \eta_{2,j} \, ,
\end{align}
where $\mu_1, \mu_2 \in C([0,1])$ and  $(\eta_{1,j})_{j\in \N}$, $(\eta_{2,j})_{j\in \N}$ are centred $C([0,1])$-valued processes satisfying the following assumption.

\begin{assumption} \label{assumption-TSP}
  The processes $(\eta_{1,j})_{j\in \N}$, $(\eta_{2,j})_{j\in \N}$ are independent centred strictly stationary processes satisfying
   Assumption~\ref{assumption} with metric  $\rho (s,t) = |t-s|^{\theta}$ for some $\theta >0 $ such that $J\theta >2$.
\end{assumption}

In the following let
\begin{align*}
C_1(s,t) &= \E[\eta_{1,j}(s) \eta_{1,j}(t)] = \cov(X_1(s),X_1(t)) \, , \\
C_2(s,t) &= \E[\eta_{2,j}(s) \eta_{2,j}(t)] = \cov(Y_1(s),Y_1(t))
\end{align*}
denote the covariance operator of the first and the second sample, respectively. We
measure the difference between $C_1$ and $C_2$ by their maximal deviation
\begin{align} \label{dev}
 d_\infty =  \| C_1 - C_2 ||_\infty=
\sup_{s,t \in [0,1]} |C_1(s,t) - C_2(s,t) | \, ,
\end{align}
and
are interested in testing if there exists a relevant difference between the
covariance operators, that is
\begin{align} \label{eq:relevant-hypothesis}
H^\Delta_0: d_\infty
\  \leq \Delta
\qquad \mbox{versus} \qquad H^\Delta_1: d_\infty > \Delta \, ,
\end{align}
where $\Delta\in\R$ is a pre-specified constant. Note that the classical
hypotheses
\begin{align} \label{class}
H_0: C_1 = C_2 \qquad \mbox{versus} \qquad H_1: C_1 \neq C_2
\end{align}
are obtained for the choice $\Delta = 0$.

We denote by
$\tilde X_{m,j} = X_j - \bar X_m, \, \tilde Y_{n,i} = Y_i - \bar Y_n$ the
centred random curves (here $\bar X_m$ and $\bar Y_n$ denote the mean
in the first and second sample, respecively), and estimate  the maximal deviation  $d_\infty$ in \eqref{dev} between the
two covariance operators  by
\begin{align} \label{eq:d-hat}
  \hat d_\infty:= \sup_{s,t \in [0,1]}
  \bigg|\frac{1}{m-1} \sum_{j=1}^m \tilde X_{m,j}^{\check \otimes 2} (s,t) -
  \, \frac{1}{n-1} \sum_{j=1}^n \tilde Y_{n,j}^{\check \otimes 2}(s,t) \bigg| \, .
\end{align}
Now a  reasonable decision rule is to reject the null hypothesis in
\eqref{eq:relevant-hypothesis} or \eqref{class} for large values of
$\hat d_\infty$. Our first result provides the asymptotic properties of
the statistic $ \hat d_\infty$.

\begin{prop}
\label{prop1}
If $\mu_1, \mu_2 \in C([0,1])$ and  $(\eta_{1,j})_{j\in \N}$, $(\eta_{2,j})_{j\in \N}$ are strictly stationary
and centred $C([0,1])$-valued processes satisfying Assumption  \ref{assumption-TSP}
and $\frac {m}{m+n} \longrightarrow \lambda \in (0,1)$ as $m,n \to \infty$, the following assertions hold true.
 \begin{itemize}
\item[(1)] If $d_\infty =0$, then
\begin{align} \label{eq:sup_T}
  \sqrt{m+n} \, \hat d_\infty \stackrel{\mathcal{D}}{\longrightarrow}
  T = \sup_{s,t \in [0,1]} |Z(s,t)| \, ,
\end{align}
 where  $Z$ is a  Gaussian random element in $C([0,1]^2)$
with covariance operator
\begin{align} \label{eq:Z-cov}
\mathbb{C} = \frac{1}{\lambda} ~ \mathbb{C}_1
-  \frac{1}{1-\lambda} ~ \mathbb{C}_2 \,  ,
\end{align}
and
$\mathbb{C}_1$ and $\mathbb{C}_2$ are the  long-run covariance operators defined by
\begin{align} \label{c1}
  \mathbb{C}_1((s,t),(s^\prime, t^\prime)) &= \sum_{i=-\infty}^{\infty} \cov\big(\eta_{1,0}^{\check \otimes 2}(s,t),
  \eta_{1,i}^{\check \otimes 2}(s^\prime,t^\prime)\big) \, , \\
  \label{c2}
  \mathbb{C}_2((s,t),(s^\prime, t^\prime)) &= \sum_{i=-\infty}^{\infty} \cov\big(\eta_{2,0}^{\check \otimes 2}(s,t),
  \eta_{2,i}^{\check \otimes 2}(s^\prime,t^\prime)\big) \, .
\end{align}
\item[(2)] If $d_\infty  >0$, we have
\begin{align} \label{eq:T-limit}
  \sqrt{m+n} \, (\hat d_\infty - d_\infty)
  \stackrel{\mathcal{D}}{\longrightarrow} T({ \cal E})
  = \max \Big\{ \sup_{(s,t) \in { \cal E}^+} Z(s,t),
  \sup_{(s,t) \in  { \cal E}^-} - Z(s,t) \Big\} \, ,
\end{align}
  where  $Z$ is a  Gaussian random element in $C([0,1]^2)$
with covariance operator defined by \eqref{eq:Z-cov} and
\begin{align} \label{eq:extremal-sets}
  \mathcal{E}^{\pm} =
  \big \{ (s,t) \in [0,1]^2  \colon C_1(s,t) - C_2(s,t) = {\pm} d_\infty \big\}
\end{align}
are the extremal sets of the difference of the covariance operators $C_1, C_2$.
\end{itemize}
\end{prop}

If $u_{1 - \alpha}$  denotes the   $(1-\alpha)$-quantile of the distribution of the random variable $T$
defined in \eqref{eq:sup_T}, a  consistent and asymptotic level  $\alpha$ tests  for the classical hypotheses in \eqref{class}
can be obtained by  rejecting the null hypothesis, whenever
\begin{align*} 
  \hat d_\infty > \frac {u_{1 - \alpha}}{\sqrt{m+n}} \, .
\end{align*}
Similarly, the null hypothesis in \eqref{eq:relevant-hypothesis} is rejected if
\begin{align*} 
  \hat d_\infty > \Delta + \frac {u_{1 - \alpha, { \cal E}}}{\sqrt{m+n}}
\end{align*}
where $u_{1 - \alpha, { \cal E}}$ is the $\alpha$-quantile of the distribution of the random variable $T(\cal E)$ defined in \eqref{eq:T-limit}.
However, the quantile $u_{1 - \alpha}$  depends on the long-run covariance operators $\mathbb{C}_1$ and $\mathbb{C}_2$ which are difficult to estimate.
For the problem of testing  relevant hypotheses the situation is even more complicated as the quantile $u_{1 - \alpha, { \cal E}}$
additionally depends on the unknown extremal sets defined in \eqref{eq:extremal-sets}, which have to be estimated as well. To deal with these problems
 we propose a bootstrap approach, which is explained for the classical  and relevant hypotheses separately.

\subsection{Classical hypotheses} \label{sec31}

In order to avoid the problem of estimating the long-run covariance operators we propose a
bootstrap procedure   to mimic the covariance structure of the distribution
of the process
$$
\frac{1}{m-1} \sum_{j=1}^m \tilde X_{m,j}^{\check \otimes 2}  -
  \, \frac{1}{n-1} \sum_{j=1}^n \tilde Y_{n,j}^{\check \otimes 2} - (C_1 -C_2)
$$
by  a multiplier bootstrap process (note that the second term vanishes in the case $d_\infty=0$). To be precise, we denote by
$(\xi_k^{(1)})_{k\in\mathbb{N}},
\ldots ,$~$(\xi_k^{(R)})_{k\in\mathbb{N}}$
and
$(\zeta_k^{(1)})_{k\in\mathbb{N}},
\ldots ,$~$(\zeta_k^{(R)})_{k\in\mathbb{N}}$
independent sequences of independent standard normal distributed
random variables
 and define the
$C([0,1]^2)$-valued processes
$\hat B_{m,n}^{(1)},$ $\ldots ,$ $\hat B_{m,n}^{(R)}$ by
\begin{align} \label{eq:bootProcess}
\begin{split}
\hat B_{m,n}^{(r)} = \sqrt{n+m} \bigg\{
	&\frac{1}{m} \sum_{k=1}^{m-l_1+1} \frac{1}{\sqrt{l_1}}\bigg(
  \sum_{j=k}^{k+l_1-1} \tilde X_{m,j}^{\check \otimes 2}
	-\frac{l_1}{m}\sum_{i=1}^m \tilde X_{m,i}^{\check \otimes 2}  \bigg)
	\xi_k^{(r)} \\
- &\frac{1}{n} \sum_{k=1}^{n-l_2+1} \frac{1}{\sqrt{l_2}}\bigg(
  \sum_{j=k}^{k+l_2-1} \tilde Y_{n,j}^{\check \otimes 2}
	-\frac{l_2}{n}\sum_{i=1}^n \tilde Y_{n,i}^{\check \otimes 2} \bigg)
	\zeta_k^{(r)} \bigg\}  ~~~~(r=1, \ldots , R) \, .
\end{split}
\end{align}
The parameters $l_1,l_2\in\mathbb{N}$ define  the block length such that $l_1/m\to 0$ and $l_2/n\to 0$ as $l_1,l_2,m,n\to\infty$.
Note that the dependence on $l_1$ and $l_2$ is not reflected in the notation of
the bootstrap processes. With these notations we
define the bootstrap statistics
\begin{align} \label{eq:bootStat-classic}
  T_{m,n}^{(r)} =  \sup_{s,t \in [0,1]} |\hat{B}_{m,n}^{(r)}(s,t) |  ~~~~(r=1,\dots,R) \, ,
\end{align}
and denote by $T_{m,n}^{\{\lfloor R(1-\alpha)\rfloor\}}$  the empirical $(1-\alpha)$-quantile of the bootstrap sample
$T_{m,n}^{(1)},\dots,T_{m,n}^{(R)}$. Then, rejecting the classical null hypothesis of equal covariance operators whenever
\begin{align} \label{eq:classical-boottest}
  \hat{d}_\infty > \frac{T_{m,n}^{\{\lfloor R(1-\alpha)\rfloor\}}}{\sqrt{m+n}}
\end{align}
defines a  bootstrap test for the classical hypotheses in  \eqref{class}. The following result provides  the  statistical properties of this test.

\begin{theorem}\label{thrm:bootstrap_class}
Suppose that  the error processes $(\eta_{1,j})_{j\in\mathbb{N}}$ and $(\eta_{2,j})_{j\in\mathbb{N}}$  in the representation \eqref{mod}
satisfy Assumption~\ref{assumption-TSP}.
 Let $\hat B_{m,n}^{(1)}, \ldots , \hat B_{m,n}^{(R)}$
denote the bootstrap processes defined by \eqref{eq:bootProcess} such that
$l_1 = m^{\beta_1}$, $l_2 = n^{\beta_2}$ with
\begin{align*}
  0 < \beta_i < \nu_i / (2+\nu_i)~,~~
  \bar{\tau}_i > (\beta_{i}( 2 + \nu_{i}) + 1) / (2 + 2\nu_{i})
\end{align*}
where $\bar \tau_i, \nu_i$ are given in Assumption \ref{assumption}, $i=1,2$.

Then, under the classical null hypothesis $H_{0}:C_{1}=C_{2}$ in \eqref{class}
we have
\begin{align} \label{eq:boot_alpha_classic}
  \lim_{m,n,R\to\infty} \mathbb{P}\bigg ( \hat{d}_{\infty}
  >  \frac{T_{m,n}^{\{\lfloor R(1-\alpha)\rfloor\}}}{\sqrt{m+n}} \bigg)
  = \alpha \, .
\end{align}
Under the alternative $H_{1}:C_{1} \not =C_{2}$ in \eqref{class}   it follows for any $R \in \N$,
\begin{align} \label{eq:consistency}
  \liminf_{m,n\to\infty} \mathbb{P}\bigg( \hat{d}_{\infty}
  >  \frac{T_{m,n}^{\{\lfloor R(1-\alpha)\rfloor\}}}{\sqrt{m+n}} \bigg)
  =1 \, .
\end{align}
\end{theorem}

\subsection{Relevant hypotheses} \label{sec22}

For testing relevant hypotheses it is crucial to estimate
the extremal sets in \eqref{eq:extremal-sets} properly. For this purpose we propose
\begin{align} \label{eq:estimatedSets}
  \hat{\mathcal{E}}_{m,n}^\pm
  &= \Big\{(s,t)\in [0,1]^2 \colon \pm \big(\hat{C}_1(s,t)-\hat{C}_2(s,t) \big)
  \geq \hat{d}_\infty - \frac{c_{m,n}}{\sqrt{m+n}} \ \Big\}
\end{align}
as estimators of the sets $\cal E^\pm$ where
$(c_{m,n})_{m,n\in\N}$ is a sequence of positive constants satisfying $\lim_{m,n \to \infty} c_{m,n} / \log(m+n) = c$
for some $c >0$. For the construction of a test of the relevant hypotheses in
\eqref{eq:relevant-hypothesis} we recall the definition of the bootstrap process in \eqref{eq:bootProcess}
and define the  statistics
\begin{align} \label{eq:bootStat}
K^{(r)}_{m,n}
=  \max\Big\{ \sup_{(s,t)\in \hat{\mathcal{E}}_{m,n}^+}
\hat{B}_{m,n}^{(r)}(s,t), \ \sup_{(s,t)\in \hat{\mathcal{E}}_{m,n}^-}
\big(- \hat{B}_{m,n}^{(r)}(s,t) \big) \Big\}  \quad \quad (r=1,\dots,R) \,
\end{align}
which serves as the bootstrap analogue of the  statistic $T({ \cal E})$ defined in
\eqref{eq:T-limit}.
  If $K_{m,n}^{\{\lfloor R(1-\alpha)\rfloor\}}$ denotes  the empirical
$(1-\alpha)$-quantile of the bootstrap sample
$K_{m,n}^{(1)},\dots,K_{m,n}^{(R)}$ we propose to reject  the null hypothesis of no relevant
difference in the covariance operators  at level $\alpha$,  whenever
\begin{align} \label{eq:relevant-boottest}
  \hat{d}_{\infty} > \Delta
  + \frac{K_{m,n}^{\{\lfloor R(1-\alpha)\rfloor\}}}{\sqrt{m+n}} \, .
\end{align}
The final result of this section states that this test
is consistent and has asymptotic level $\alpha$.

\begin{theorem}\label{thrm:bootstrap_test}
Suppose that the assumptions of Theorem \ref{thrm:bootstrap_class} are satisfied and that $\Delta >0$.
\begin{itemize}
\item[(1)]
Under the null hypothesis $H_0: d_\infty \leq \Delta$ of no relevant difference in the covariance  operators, it follows
\begin{align} \label{eq:boot_alpha}
  \lim_{m,n,R\to\infty} \mathbb{P}\bigg ( \hat{d}_{\infty}
  > \Delta + \frac{K_{m,n}^{\{\lfloor R(1-\alpha)\rfloor\}}}{\sqrt{m+n}} \bigg)
  = \alpha \, ,
\end{align}
  if $\Delta = d_\infty$ and, for any $R \in \N$,
\begin{align*}
  \lim_{m,n\to\infty} \mathbb{P}\bigg ( \hat{d}_{\infty}
  > \Delta + \frac{K_{m,n}^{\{\lfloor R(1-\alpha)\rfloor\}}}{\sqrt{m+n}} \bigg)
  = 0 \, ,
\end{align*}
if $0 < d_\infty < \Delta$.

\item[(2)]Under the alternative $H_1: d_\infty > \Delta$  of a relevant difference in
the covariance operators it follows for any $R \in \N$
\begin{align*}
  \liminf_{m,n\to\infty} \mathbb{P}\bigg( \hat{d}_{\infty}
  > \Delta + \frac{K_{m,n}^{\{\lfloor R(1-\alpha)\rfloor\}}}{\sqrt{m+n}} \bigg)
  =1 \, .
\end{align*}
\end{itemize}
\end{theorem}

\section{Detecting changes in the covariance operator} \label{sec4}
\def\theequation{4.\arabic{equation}}
\setcounter{equation}{0}

In this section we study the change point problem for the covariance operator
of an array $(X_{n,j}\colon n\in \N,j=1,\dots,n)$ of $C([0,1])$-valued random variables.
For the consideration of relevant changes we require a dependence concept for an
array $(\tilde \eta_{n,j} \colon n\in\mathbb{N}, j = 1,\dots, n)$ of
random variables in $C(T)$ with strictly stationary rows. For this purpose we   denote by $\mathcal{F}^{k^\prime}_{k,n}$ the
$\sigma$-field generated by $(\tilde \eta_{n,j}\colon k\leq j \leq k^\prime)$.
The $k$th $\varphi$\textit{-mixing coefficient} of the array
$(\tilde \eta_{n,j} \colon n\in\mathbb{N}, j = 1,\dots, n)$ is then defined by
\[
\varphi (k)
= \sup_{n \in \mathbb{N}} \sup_{k^\prime \in \{1,\dots,n-k\}} \phi (\mathcal{F}_{1,n}^{k^\prime} ,\mathcal{F}_{k^\prime+k,n}^n)
\]
and $(\tilde \eta_{n,j} \colon n\in\mathbb{N}, j = 1,\dots, n)$ is
called $\varphi$\textit{-mixing} whenever $\varphi(k) \to 0$ as $k\to\infty$.
For our theoretical investigations we make the following assumption.

\begin{assumption} \label{assumption-CPP}
For some $\vartheta\in (0, \frac{1}{2}]$ there exists a number $s^* \in [\vartheta,1-\vartheta]$ such
that  the random variables $(X_{n,j}\colon n\in \N,j=1,\dots,n)$ are given by $X_{n,j} = \mu + \tilde{\eta}_{n,j}$, where $\mu = \E[X_{n,j}]$ denotes the common expectation function,
\begin{align} \label{eq:tilde-eta}
  \tilde \eta_{n,j} =
  \begin{cases}
    \eta_{1,j} \quad \text{ if } ~j\in \{1,\dots,\lfloor s^* n \rfloor\} \\
    \eta_{2,j} \quad \text{ if } ~j\in \{\lfloor s^* n \rfloor+1,\dots,n \}
  \end{cases}
\end{align}
and $(\eta_{1,j})_{n\in\N}$, $(\eta_{2,j})_{n\in\N}$ are centred strictly stationary processes satisfying conditions  (A1) - (A3) of Assumption~\ref{assumption} with metric $\rho(s,t) = |s-t|^\theta$ for some $\theta >0$ such that $\theta J > 2$. Furthermore it is assumed that the array $(\tilde \eta_{n,j} \colon n\in\mathbb{N}, j = 1,\dots, n)$ is $\varphi$\textit{-mixing} with mixing coefficients satisfying condition   (A4) of Assumption~\ref{assumption}.
\end{assumption}

We denote by $C_1$  and $C_2$ the covariance operator before and after the change point. Recalling  the definition of $d_\infty$ in \eqref{dev} the relevant and classical hypotheses
are given by \eqref{eq:relevant-hypothesis} and \eqref{class}, respectively.
 For the construction of a test for these hypotheses we consider a
 sequential empirical process on $C([0,1]^3)$  defined by
\begin{align}   \label{un}
  \hat {\mathbb{U}}_n(s,t,u) = \frac {1}{n}
  \Big( \sum^{\lfloor sn \rfloor}_{j=1}
  \tilde X_{n,j}^{\check \otimes 2}(t,u)
  + n\Big(s- \frac {\lfloor sn \rfloor}{n}\Big)
  \tilde X_{n, \lfloor sn \rfloor+1}^{\check \otimes 2}(t,u)
  - s \sum^n_{j=1} \tilde X_{n,j}^{\check \otimes 2}(t,u) \Big)
\end{align}
where $\tilde X_{n,j} = X_{n,j} - \bar X_n$ $(j = 1,\dots,n; ~ n\in\N)$ and note that it can be shown that
\begin{align*}
  \mathbb{E}  \big[ \hat {\mathbb{U}}_n (s,t,u) \big]
  = \big( s \wedge s^* - s s^*  \big) \big (C_1 (t,u) - C_2(t,u) \big)  + o_\mathbb{P}(1) \, .
\end{align*}
Consequently, it is reasonable to consider the statistic
\begin{equation} \label{eq:M-hat}
\mathbb{\hat  M}_n  =  \sup_{s \in [0,1]} \sup_{t,u \in [0,1]} | \hat {\mathbb{U}}_n (s,t,u) |
\end{equation}
as an estimate of
$$
s^*(1-s^*) \, d_\infty = s^*(1-s^*) \, \|C_1-C_2\|_\infty \, .
$$
The following result makes these heuristic arguments  precise.

\begin{prop}
\label{prop2}
 If  Assumption  \ref{assumption-CPP} is satisfied, the following statements hold true.
 \begin{itemize}
\item[(1)] If $d_\infty =0$, then
\begin{align} \label{teil1}
  \sqrt{n} \, \mathbb{\hat  M}_n \stackrel{\cal D}{\longrightarrow} \check T
  =  \sup_{(s,t,u) \in [0,1]^3} | \mathbb{W}(s,t,u)|
\end{align}
where $\mathbb{W}$ is a  Gaussian random element in $C([0,1]^3)$
with covariance operator
\begin{align} \label{eq:W-cov}
\begin{split}
\cov& (\mathbb{W}(s,t,u), \mathbb{W}(s^\prime,t^\prime, u^\prime)) \\
=& \big\{(s\wedge s^\prime \wedge s^*) + s s^\prime s^* - s(s^\prime \wedge s^*) - s^\prime(s \wedge s^*) \big\} \, \mathbb C_1((t,u),(t^\prime,u^\prime)) \, \\
&+ \big\{(s\wedge s^\prime - s^*)_+ + s s^\prime(1-s^*) - s(s^\prime - s^*)_+ - s^\prime (s - s^*)_+ \big \} \, \mathbb C_2((t,u),(t^\prime,u^\prime))
\end{split}
\end{align}
and the long-run covariance operators $\mathbb{C}_1, \mathbb{C}_2$ are defined by
\begin{align} \label{eq:lrv-cp}
\mathbb{C}_l((s,t),(s^\prime, t^\prime)) &= \sum_{i=-\infty}^{\infty} \cov\big(\eta_{l,0}^{\check \otimes 2}(s,t),
\eta_{l,i}^{\check \otimes 2}(s^\prime,t^\prime)\big) \quad \quad (l=1,2) \, .
\end{align}
\item[(2)] If $d_\infty  >0$, we have
\begin{align} \label{eq:D-limit}
\begin{split}
   \sqrt{n} \big( \mathbb{\hat  M}_n
  - s^*(1-s^*) d_\infty \big)
  \stackrel{\cal D}{\longrightarrow}
  \tilde D (\mathcal{E})
  = \max \Big \{ \sup_{(t,u) \in \mathcal{E}  ^+} \mathbb{W}(s^*,t,u),
  \sup_{(t,u) \in \mathcal{E} ^-} - \mathbb{W}(s^*,t,u) \Big \} \, ,
\end{split}
\end{align}
 where   $\mathbb{W}$ is a  Gaussian random element in $C([0,1]^3)$
with covariance operator defined by \eqref{eq:W-cov}
and  $  \mathcal{E}^{\pm}$ are the extremal sets
defined in \eqref{eq:extremal-sets}.
\end{itemize}
\end{prop}

As in the two sample problem we can form decision rules, rejecting  the null hypothesis (classical or relevant) for large values
of  $\mathbb{\hat  M}_n$. Note that this requires estimation of the long-run covariance operators and (in the case of relevant hypotheses)
the  estimation of the change point and the extremal sets. For the construction of an explicit test (based on a multiplier bootstrap) we investigate
 again classical  and relevant hypotheses separately.

\subsection{Classical hypotheses} \label{sec41}

Most of the literature on change point analysis of covariance operators investigates the classical hypotheses of
the form  \eqref{class}, where  $C_1$ and $C_2$ denote the covariance operator before and after the change point
 \citep[see][]{jaruskova2013, sharipov2019, stoehr2019}. In order to obtain critical values for a test for a structural break in the covariance operators
 we consider a
$C([0,1]^3)$-valued bootstrap process defined by
\begin{align} \label{eq:bootProcess-cp}
\begin{split}
  \hat{B}_n^{(r)}(s,t,u) =& \frac{1}{\sqrt{n}} \sum_{k=1}^{\lfloor sn \rfloor}
  \frac{1}{\sqrt{l}} \Big( \sum_{j=k}^{k+l-1} \hat{Y}_{n,j}(t,u)
  - \frac{l}{n} \sum_{j=1}^n \hat{Y}_{n,j}(t,u) \Big) \xi_k^{(r)} \\
  &+ \sqrt{n}\Big(s - \frac{\lfloor sn \rfloor}{n} \Big)\frac{1}{\sqrt{l}}
  \Big( \sum_{j=\lfloor sn \rfloor +1}^{\lfloor sn \rfloor+l} \hat{Y}_{n,j}(t,u)
  - \frac{l}{n} \sum_{j=1}^n \hat{Y}_{n,j}(t,u) \Big)
  \xi_{\lfloor sn \rfloor +1}^{(r)} \, ,
\end{split}
\end{align}
if $\lfloor sn \rfloor \leq  n-l$, where $(\xi_k^{(1)})_{k\in\N}, \ldots, (\xi_k^{(R)})_{k\in\N}$ denote independent sequences of independent Gaussian random variables with mean $0$ and variance $1$ and
\[
\hat{Y}_{n,j} = \tilde X_{n,j}^{\check \otimes 2}(t,u) - (\hat{C}_2 - \hat{C}_1)
\mathds{1}\{j > \lfloor \hat{s}n \rfloor \} \quad \quad  (j=1,\dots,n) \, .
\]
The expressions
\begin{align*}
\hat{C}_1 = \frac{1}{\lfloor \hat{s} n\rfloor}
\sum_{j=1}^{\lfloor \hat{s}n \rfloor} \tilde X_{n,j}^{\check \otimes 2}(t,u)
\qquad \text{and} \qquad
\hat{C}_2 = \frac{1}{\lfloor (1-\hat{s}) n\rfloor}
\sum_{j=\lfloor \hat{s}n \rfloor +1}^{n} \tilde X_{n,j}^{\check \otimes 2}(t,u)
\end{align*}
are estimators of the covariance operator before and after the change point and
\begin{align} \label{cpEstimator}
\hat  s  = \max \Big \{ \vartheta, \, \min \Big\{ \frac{1}{n} \arg\max_{1\leq k <n} \big\| \hat {\mathbb{U}}_n (k/n,\cdot, \cdot) \big\|_\infty, \, 1- \vartheta \Big\} \Big\}
\end{align}
is an estimator of the unknown change location $s^*$ (note that $s^* \in (\vartheta, 1-\vartheta)$ by assumption).
In \eqref{eq:bootProcess-cp} the parameter   $l\in\mathbb{N}$ denotes the block length satisfying  $l/n\to 0$ as $l,n\to\infty$ and
for any $t,u\in[0,1]$ and any $s\in[0,1]$ such that $\lfloor sn \rfloor > n-l$ we define
$$
\hat{B}_n^{(r)}((n-l)/n,t,u) = \hat{B}_n^{(r)}(s,t,u) \, .
$$
Finally, a  bootstrap process is defined by
\begin{align} \label{eq:What-boot}
\hat{\mathbb{W}}_n^{(r)}(s,t,u)
= \hat{B}_n^{(r)}(s,t,u)-s\hat{B}_n^{(r)}(1,t,u)  \quad \quad (r = 1,\dots,R)
\end{align}
and we consider the bootstrap statistic
\begin{align} \label{eq:checkT}
\check T_n^{(r) }
= \sup_{s,t,u \in [0,1]}  \big | \hat{\mathbb{W}}_n^{(r)}(s,t,u) \big | \quad \quad (r = 1,\dots,R) \, .
\end{align}
If  $\check T_n^{\{\lfloor R(1-\alpha)\rfloor\}}$  denotes
the empirical $(1-\alpha)$-quantile of the bootstrap sample
$\check T_n^{(1)},\check T_n^{(2)} , \dots ,\check  T_n^{(R)}$, the
classical null hypothesis \eqref{class} of no change in the covariance operators is rejected, whenever
\begin{align} \label{eq:classical-test-cp}
\hat{\mathbb{M}}_n > \frac{\check T_n^{\{\lfloor R(1-\alpha) \rfloor\}}}{\sqrt{n}} \, .
\end{align}

\begin{theorem}\label{thrm:bootstrap_test_cpclass}
Assume that the array $(X_{n,j}\colon n\in \N,j=1,\dots,n)$ satisfies
Assumption~\ref{assumption-CPP}.
Further assume that  $l = n^\beta$
for some   constant $\beta \in (1/5,2/7)$ such that the constant $\nu$ in (A2) satisfies
$\nu\geq 4$ and
$$(\beta(2+\nu)+1)/(2+2\nu) < \bar \tau < 1/2
$$
where $\bar \tau $ is defined in  (A4).

Then, under the classical  null hypothesis  $H_0: C_1=C_2$, we have
\begin{align*}
  \lim_{n,R\to\infty} \mathbb{P}\bigg ( \hat{\mathbb{M}}_n
  >   \frac{\check T_{n}^{\{\lfloor R(1-\alpha)\rfloor\}}}{\sqrt{n}} \bigg)
  = \alpha \, .
\end{align*}
Under the alternative $H_1: C_1 \neq C_2$ we have, for
any $R \in \N$,
\begin{align*}
  \liminf_{n\to\infty} \mathbb{P}\bigg( \hat{\mathbb{M}}_n
  >  \frac{\check T_{n}^{\{\lfloor R(1-\alpha)\rfloor\}}}{\sqrt{n}} \bigg)
  =1 \, .
\end{align*}
\end{theorem}

\subsection{Relevant  hypotheses} \label{sec42}

Testing  for a relevant  change   in the covariance operators  as formulated in \eqref{eq:relevant-hypothesis}  is more complicated. In particular because - as indicated in Proposition
\ref{prop2}  - it  additionally requires  the estimation of the  extremal sets.
To be precise we recall the definition of $\hat{\mathbb{M}}_n$ in \eqref{eq:M-hat} and define
\begin{align} \label{eq:statistic-cp}
\hat{d}_\infty = \frac{ \mathbb{\hat  M}_n}{\hat{s}(1-\hat{s})}
\end{align}
as an estimator of the maximal deviation of the covariance operator before and after the
change point,  and use
\begin{align} \label{eq:estimatedSets-cp}
\hat{\mathcal{E}}_{n}^\pm  &= \Big\{ (t,u)\in[0,1]^2 \colon \pm ( \hat{C}_1(t,u)-\hat{C}_2(t,u) )
\geq \hat{d}_\infty - \frac{c_{n}}{\sqrt{n}} \ \Big\} \, ,
\end{align}
as the estimator of the extremal sets, where
$(c_n)_{n\in\N}$ is a sequence of positive constants such that $\lim_{n \to \infty}  c_n  / \log(n) = c  >0$.
In order to obtain a test for the relevant hypotheses  in \eqref{eq:relevant-hypothesis} define, for $r=1,\dots,R$, the bootstrap statistics
\begin{align} \label{eq:bootStat-cp}
  \check K^{(r)}_{n} =  \frac{1}{\hat{s}(1-\hat{s})}\max\big\{
  \sup_{(t,u)\in \hat{\mathcal{E}}_{n}^+} \hat{W}_{n}^{(r)}(\hat{s},t,u), \
  \sup_{(t,u)\in \hat{\mathcal{E}}_{n}^-} \big(- \hat{W}_{n}^{(r)}(\hat{s},t,u)
  \big) \big\} \, .
\end{align}
Then the null hypothesis of no relevant change in the covariance operators is  rejected at level $\alpha$, whenever
\begin{align} \label{eq:relevant-test-cp}
\hat{d}_{\infty} > \Delta + \frac{\check K_n^{\{\lfloor R(1-\alpha)\rfloor\}}}{\sqrt{n}} \, ,
\end{align}
where $\check K_n^{\{\lfloor R(1-\alpha)\rfloor\}}$ is the empirical $(1-\alpha)$-quantile of the bootstrap sample
$\check K_n^{(1)}, $ $\check K_n^{(2)} ,$ $ \dots ,\check K_n^{(R)}$.
The following result  shows that the bootstrap test for the relevant hypotheses is
consistent and has asymptotic level $\alpha$.

\begin{theorem}\label{thrm:bootstrap_test_cp}
Let  the assumption of Theorem \ref{thrm:bootstrap_test_cpclass} be satisfied and furthermore assume that the random variable $M$ in (A3) is bounded.
  \begin{itemize}
\item[(1)]
Under the null hypothesis $H_0: d_\infty \leq \Delta$ of no relevant difference in the covariance operators, we have
\begin{align*}
  \lim_{n,R\to\infty} \mathbb{P}\bigg ( \hat{d}_{\infty}
  > \Delta + \frac{\check K_{n}^{\{\lfloor R(1-\alpha)\rfloor\}}}{\sqrt{n}} \bigg)
  = \alpha \, ,
\end{align*}
if $\Delta = d_\infty$ and, for any $R \in \N$,
\begin{align*}
  \lim_{n\to\infty} \mathbb{P}\bigg ( \hat{d}_{\infty}
  > \Delta + \frac{\check K_{n}^{\{\lfloor R(1-\alpha)\rfloor\}}}{\sqrt{n}} \bigg)
  = 0 \, ,
\end{align*}
if $0 < d_\infty < \Delta$.
  \item[(2)]
Under the alternative $H_1: d_\infty > \Delta$ of a relevant difference in the covariance operators,  we have for
any $R \in \N$,
\begin{align*}
  \liminf_{n\to\infty} \mathbb{P}\bigg( \hat{d}_{\infty}
  > \Delta + \frac{\check K_{n}^{\{\lfloor R(1-\alpha)\rfloor\}}}{\sqrt{n}} \bigg)
  =1 \, .
\end{align*}
\end{itemize}
\end{theorem}

\section{Finite sample properties} \label{sec:simulations}
\def\theequation{5.\arabic{equation}}
\setcounter{equation}{0}

\subsection{Simulation study} \label{sec51}

In this section we study the finite sample properties of the test procedures
developed in this paper and
we also compare it with some  competing procedures from the literature, which
can be used under similar assumptions as considered here.  The   empirical rejection probabilities of the
different tests have been calculated by $1000$ simulation runs   and $200$ bootstrap
statistics are used for the calculation of the bootstrap quantiles in each run.

\subsubsection{Two sample problem}

\paragraph{Classical hypotheses:}
In the following we investigate the finite sample properties of the test
\eqref{eq:classical-boottest} for   the classical null hypothesis of
equal covariance operators in \eqref{class}. For the sake of comparison, we use the same
scenarios as considered in \cite{Paparoditis2016} who   developed a
  bootstrap  test for the hypotheses \eqref{class}. \cite{Paparoditis2016} also applied the FPC test developed by \cite{fremdt2013} to these scenarios, such that a comparison with the method developed by these authors is also possible.  To be  precise,  curves are
generated according to the model
\begin{align} \label{eq:non-gaussian}
\begin{split}
X_{i}(t) &= \sum_{k=1}^{10}
\big \{ 2^{1/2} k^{-1/2} \sin(\pi k t) V_{i,k}
+ k^{-1/2} \cos(2 \pi k t) W_{i,k} \big \} \\
Y_{j}(t) &= c \sum_{k=1}^{10}
\big \{ 2^{1/2} k^{-1/2} \sin(\pi k t) \tilde{V}_{j,k}
+ k^{-1/2} \cos(2 \pi k t) \tilde{W}_{j,k} \big \}
\end{split}
\end{align}
($i=1, \ldots , m, \, j=1, \ldots , n$), where the random variables
$V_{i,k}, W_{i,k}, \tilde{V}_{j,k}, \tilde{W}_{j,k}$
 are independent and
  $t_5$-distributed. The constant $c$ determines if the null hypothesis $(c=1)$ holds or not $(c \neq 1)$.
     In order to obtain functional data
objects, the curves are evaluated at $500$ equidistant points in $[0,1]$ and
then the Fourier basis consisting of $49$ basis functions is used to transform
these function values into a functional data object
(using the function ``Data2fd" from the ``fda" R-package).

In Table~\ref{tab:paparoditis} we display empirical rejection probabilities for
two different sample sizes and different choices of $c$. \cite{Paparoditis2016}
state that the procedure proposed by \cite{fremdt2013} achieves the
best results  if two FPCs are used to represent the data, and therefore,
the results of this procedure were obtained for this case.

\begin{table}[t!]
  \begin{center}
    \setlength{\tabcolsep}{5pt}
    { \scriptsize
    \begin{tabular}{c|ccc||ccc}
      & \multicolumn{3}{|c||} {$c = 1$} &    \multicolumn{3}{|c} {$c = 1.2$} \\
      \hline
      $n,m$                & 1\%      & 5\%        & 10\%       &  1\%          & 5\%          & 10\%  \\
      \hline \hline
      \multirow{2}{*}{25}  & 0.9      & 4.2        & 11.8       &   3.0       & 13.4       & 24.7        \\
                           & (0, 0.3) & (0.6, 2.5) & (2.2, 8.2) &  (0, 0.5)   & (1.6, 5.0) & (3.9, 14.7) \\[4pt]
      \multirow{2}{*}{50}  & 0.8      & 3.6        & 8.6        &  6.6        & 22.4       & 35.0        \\
                           & (0, 0.6) & (1.6, 3.2) & (4.1, 7.6) &  (0.3, 0.8) & (2.6, 9.8) & (7.2, 23.9) \\[4pt]
      \hline \hline
      & \multicolumn{3}{|c||} {$c = 1.4$} &    \multicolumn{3}{|c} {$c = 1.6$} \\
      \hline
      $n,m$                & 1\%         & 5\%         & 10\%         &  1\%         & 5\%          & 10\%        \\
      \hline \hline
      \multirow{2}{*}{25}  & 10.3        & 32.3        & 51.0         &  22.4        & 54.7         & 73.8         \\
                           & (0, 1.6)    & (1.1, 16.8) & (5.2, 36.8)  &  (0, 4.7)    & (1.0, 33.8)  & (9.5, 61.2)  \\[4pt]
      \multirow{2}{*}{50}  & 27.8        & 58.9        & 75.1         &  55.0        & 83.3         & 91.8         \\
                           & (0.2, 12.8) & (6.5, 46.1) & (22.1, 67.6) &  (1.4, 37.0) & (28.5, 79.6) & (55.9, 90.3) \\[4pt]
      \hline \hline
      & \multicolumn{3}{|c||} {$c = 1.8$} &    \multicolumn{3}{|c} {$c = 2$} \\
      \hline
      $n,m$                 & 1\%        & 5\%          & 10\%         &  1\%          & 5\%          & 10\%        \\
      \hline \hline
      \multirow{2}{*}{25}  & 34.9        & 72.2         & 87.4         &  45.3         & 81.9         & 93.3         \\
                           & (0, 10.4)   & (3.6, 55.7)  & (23.0, 82.3) &  (0, 17.7)    & (7.0, 66.6)  & (50.5, 89.2)  \\[4pt]
      \multirow{2}{*}{50}  & 73.0        & 93.4         & 97.5         &  83.0         & 96.4         & 98.6         \\
                           & (6.6, 61.2) & (57.4, 91.5) & (82.1, 96.6) &  (24.5, 74.2) & (83.6, 93.7) & (95.7, 97.7) \\[4pt]
    \end{tabular}
    }
    \caption{\it Rejection probabilities of the test
      \eqref{eq:classical-boottest}  for the classical hypotheses \eqref{class}. The case $c=1$
      corresponds to the null hypothesis.
      The numbers in the brackets display the empirical rejection probabilities
      of the tests proposed by \cite{fremdt2013} and \cite{Paparoditis2016},
      respectively.} \label{tab:paparoditis}
  \end{center}
\end{table}

We observe that under the null, i.e. $c=1$, the nominal level is
 well approximated by the test \eqref{eq:classical-boottest} and the alternatives are detected with reasonable probability. Moreover,   in
all considered scenarios under the alternative, the new procedure achieves a better
power than the tests of \cite{Paparoditis2016} and \cite{fremdt2013}.

\paragraph{Relevant hypotheses:} \label{sec:simulations-relevant}


We now investigate
the finite sample properties of the decision rule \eqref{eq:relevant-boottest} for testing relevant hypotheses of the form
\eqref{eq:relevant-hypothesis} in the two sample problem. For this purpose we define different processes
including independent random functions, functional moving average processes and
non-Gaussian random curves.

For the data generation, we proceed similarly as in Sections~6.3
and 6.4 of \cite{aueDubartNorinhoHormann2015}. We consider $21$
$B$-spline basis functions $\nu_1,\dots , \nu_{21}$  and restrict
to functions in the linear space
$\mathbb{H}=\mathrm{span}\{ \nu_1,\dots , \nu_{21} \} $.
Then, for a sample of size $m\in\N$, random functions
$\varepsilon_1, \dots, \varepsilon_m \in \mathbb{H} \subset  C([0,1])$
are defined by
\begin{equation}\label{eq:errors}
  \varepsilon_j = \sum_{i=1}^{21} N_{i,j} \nu_i \, ,
  \qquad j=1,\dots , m \, ,
\end{equation}
where \(N_{1,j},N_{2,j}, \dots, N_{21,j} \) are independent normally
distributed random variables with expectation zero and variance
$\var(N_{i,j})=\sigma_i^2=1/i^2$.
Independent and identically distributed Gaussian random functions are then obtained by
\begin{align} \label{eq:fIID}
  X_i = \varepsilon_i \qquad (i=1,\ldots,m) \, ,
\end{align}
and we call $\{X_i\}^m_{i=1}$ fIID process.
In order to obtain independent non-Gaussian curves,   we replace the normally distributed random coefficients in \eqref{eq:errors} by independent $t_5$-distributed random variables, that is
 $N_{i,j} \sim t_5 \, \sqrt{3/(5i^2)}$.
Then, the variances of the coefficients are the same as for the fIID processes
and the corresponding setting is called the non-Gaussian process.

Using the processes in \eqref{eq:errors}, fMA(2) processes can be defined by
\begin{align} \label{eq:fMA2}
  X_i = \varepsilon_i + \kappa_1 \, \varepsilon_{i-1}
  + \kappa_2 \, \varepsilon_{i-2} \qquad (i=1,\ldots,m)
\end{align}
where $\kappa_1, \kappa_2 \in \mathbb{R}$ are
parameters defining the dependency (for initialization  define $\varepsilon_{-1}, \varepsilon_0$ as  independent copies of $\varepsilon_1$). In the
simulations, we set $\kappa_1 = 0.7, \kappa_2 = 0$ to obtain an fMA(1) processes and
$\kappa_1 = 0.5, \kappa_2 = 0.3$ for an  fMA(2) processes.

In order to test for a relevant difference in the covariance operators of two
populations, we generate an independent second sample, $\tilde Y_1,\dots, \tilde Y_n$, in
the same way  and multiply it by a
constant $a$ such that   $Y_i = a \, \tilde{Y}_i$   ($i = 1,\dots,n$).
Consequently,
\begin{align}\label{h00}
  |C_1(s,t) - C_2(s,t)| = |C_1(s,t) (a^2-1)|
\end{align}
where $C_1, C_2$ are the covariance operators of $X_1$ and $Y_1$, respectively.

In the case of fIID and non-Gaussian processes defined by \eqref{eq:fIID}, the maximum of the covariance
operator is given by
\begin{align*}
  \max_{s,t\in[0,1]} \cov(X_1(s),X_1(t))
  = \max_{s,t\in[0,1]} \sum_{i=1}^D \nu_i(s)\nu_i(t) / i^2 = 1
\end{align*}
which is attained at the point $(s,t)=(0,0)$.
Consequently, we obtain for the sup-norm
\begin{align*}
 \| C_1 - C_2 \|_\infty  =  |a^2 - 1|
\end{align*}
in both cases and the extremal sets are defined by $\mathcal{E}^+ = \{(0,0)\}, \mathcal{E}^- = \emptyset$. For fMA(2) processes of the form \eqref{eq:fMA2}, we obtain
\begin{align*}
  \| C_1 - C_2 \|_\infty
  = |a^2 - 1| \, (1+\kappa_1^2 + \kappa_2^2) \, .
\end{align*}
 In Table~\ref{tab:relevant-lvl} we display empirical
rejection probabilities for the hypotheses in \eqref{eq:relevant-hypothesis} for the different types of processes and different
choices of the sample sizes. In each case, we use $a = \sqrt{2}$ and define
$\Delta$ such that  $\Delta = |a^2 - 1|$ in the fIID and
non-Gaussian setting and $\Delta = |a^2 - 1| \, (1+\kappa_1^2+\kappa_2^2)$ in
the fMA(1) and fMA(2) setting. Throughout this section we call this situation the boundary of the hypotheses
 \eqref{eq:relevant-hypothesis}.
 For the estimation of the extremal sets, we use
$c_{m,n} = 0.1 \log(n + m)$ in \eqref{eq:estimatedSets} and the  block lengths in the bootstrap process \eqref{eq:bootProcess} are chosen as
  $l_1 = l_2 = 1$ in the fIID   cases, as
$l_1 = l_2 = 2$ in the fMA(1) and as $l_1 = l_2 = 3$ in the fMA(2) case.
We observe a reasonable approximation of the nominal level of the test at the boundary
of the hypotheses in all cases under consideration. The nominal level in the interior
of the hypotheses, that is $\| C_1 - C_2 \|_\infty < \Delta$ is usually much smaller (these
results are not displayed).

\begin{table}[t]
  \centering  \setlength{\tabcolsep}{5.5pt}
      { \scriptsize
  \begin{tabular}{c|ccc||ccc||ccc||ccc}
    & \multicolumn{3}{|c||} {fIID} & \multicolumn{3}{|c||} {non-Gaussian} & \multicolumn{3}{|c||} {fMA(1)} & \multicolumn{3}{|c} {fMA(2)} \\
    \hline
    \(m,n\)     & 1\% & 5\% & 10\% & 1\% & 5\% & 10\% & 1\% & 5\% & 10\% & 1\% & 5\% & 10\%  \\
    \hline \hline
    50, 50      & 1.0 & 4.6 & 8.7  & 0   & 1.7 & 5.4  & 1.2 & 5.2 & 10.5 & 1.6 & 5.1 & 10.1  \\
    50, 100     & 0.9 & 4.7 & 10.1 & 0.6 & 3.9 & 9.8  & 2.1 & 7.6 & 13.5 & 2.4 & 7.6 & 11.9  \\
    100, 100    & 0.9 & 3.9 & 9.1  & 0.5 & 3.1 & 9.0  & 1.4 & 5.7 & 10.8 & 1.0 & 4.2 & 10.8
  \end{tabular}
  }
  \caption{\it \label{tab:relevant-lvl}
    Simulated level of the test \eqref{eq:relevant-boottest} for a relevant difference in the covariance
    operators at the boundary of the hypotheses in \eqref{eq:relevant-hypothesis}, that is $\| C_1 - C_2\|_\infty=\Delta$.}
\end{table}



\medskip

\begin{figure}[h]
  {  \centering
    \includegraphics[width=6cm,height=6cm]{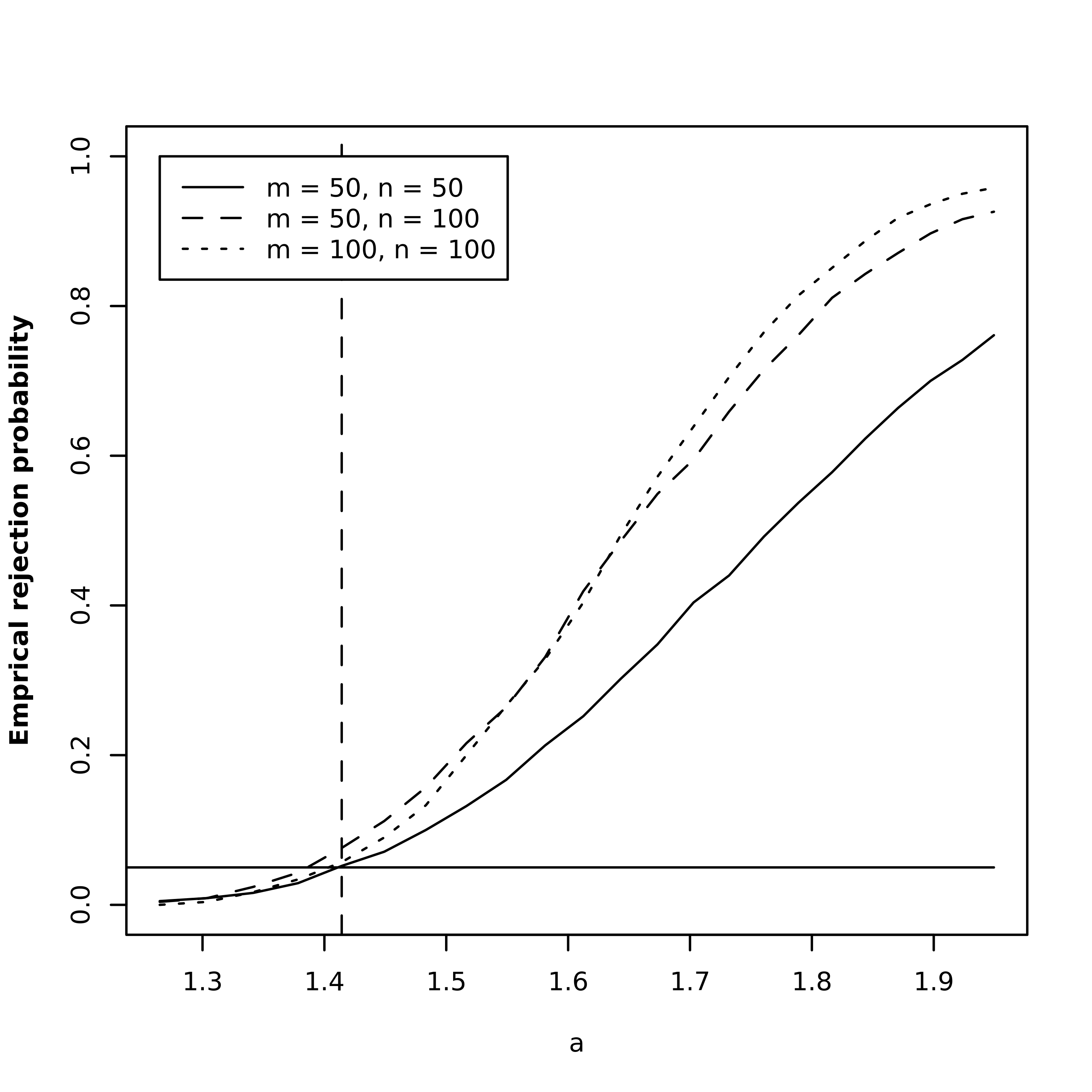}
    ~~~    ~~~     ~~~
    \includegraphics[width=6cm,height=6cm]{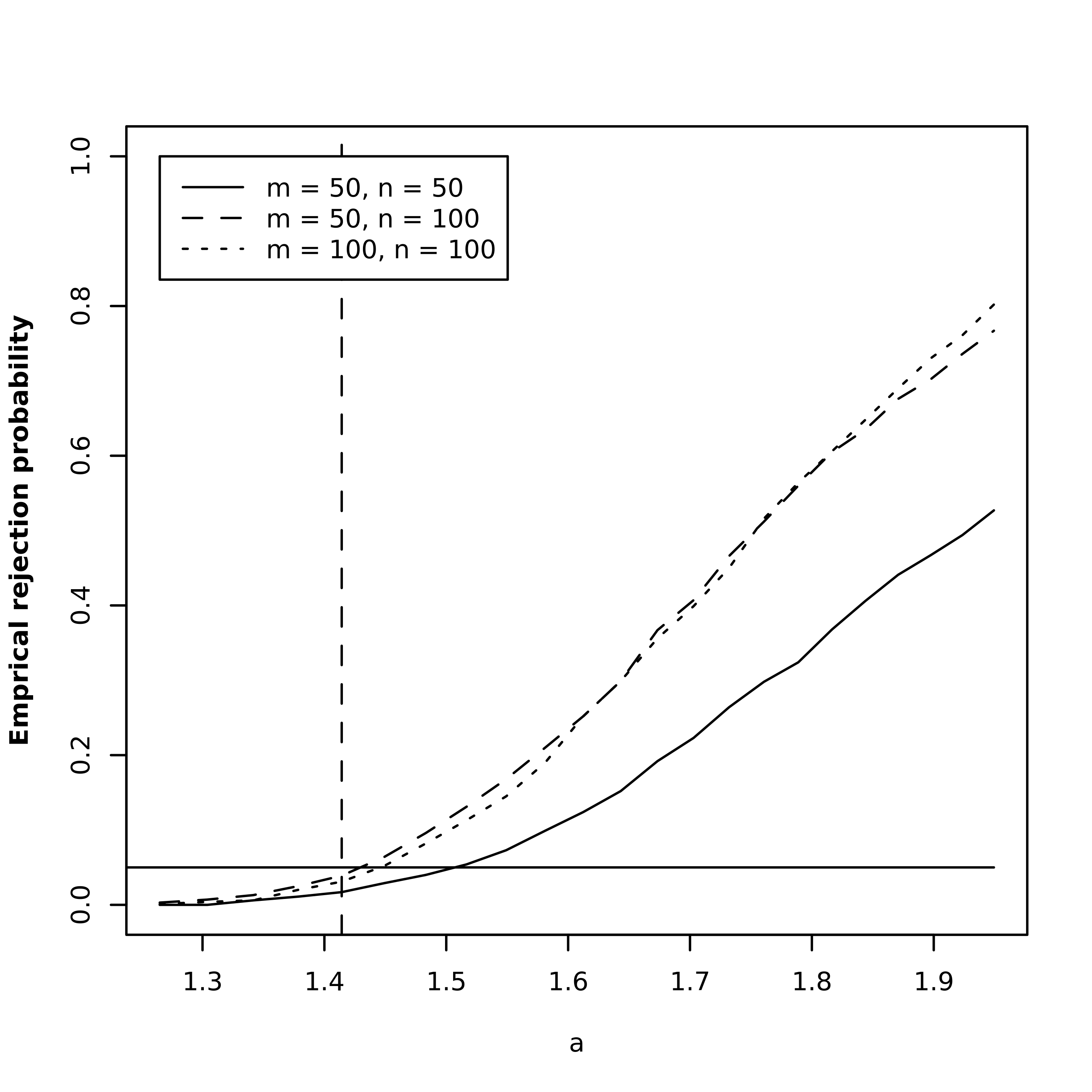}
    \caption{\label{fig:powercurves}
      \it Simulated rejection probabilities of the test \eqref{eq:relevant-boottest} for a non-relevant
      difference in the covariance operators. Left panel
        fMA(1) with threshold $\Delta=1+0.7^2$. Right panel    non-Gaussian   curves with    threshold  $\Delta = 1$. The second sample
      is multiplied by $a$ for $a = \sqrt{1.6},\sqrt{1.7},\dots,\sqrt{3.8}$, and the vertical lines represent the boundary of the null hypotheses {(i.e. $a = \sqrt{2}$)}.
    }
  }
\end{figure}

Next we study the properties of the test  \eqref{eq:relevant-boottest}
under the alternative in \eqref{eq:relevant-hypothesis}. As before  two independent
identically distributed
samples are generated  where
  the second sample is multiplied by a factor $a$. The threshold
 $\Delta$ is fixed and then empirical rejection probabilities are
simulated for different choices of the constant $a$,
  such that the properties   stated in Theorem~\ref{thrm:bootstrap_test}
can be visualized.  The results are displayed in Figure~\ref{fig:powercurves} for fMA(1)
processes (with $\kappa_1 = 0.7, \kappa_2 = 0$) and non-Gaussian random curves. The threshold  in \eqref{eq:relevant-hypothesis}   is set to $\Delta = 1+\kappa_1^2$ and
$\Delta = 1$, respectively. As illustrated before, the nominal level is reasonably
well approximated in both cases and with increasing factor $a$, the
empirical rejection probability also increases. It can be observed that the
empirical rejection probability increases slightly faster in the fMA(1) case. An explanation of this
behaviour consists in the fact
  that for the same factor $a$, the true maximal difference of
the covariance operators is greater in the fMA(1) than in the non-Gaussian
case. 

\subsubsection{Change point problem}

\paragraph{Classical hypotheses:}


We begin with a comparison of the test
 \eqref{eq:classical-test-cp} for the classical hypotheses \eqref{class} with
 two procedures which were recently proposed by \cite{sharipov2019} and are based on the sup and $L^2$-norm of the CUSUM statistic.
Following these authors we generate data from   the model
\begin{align} \label{eq:sharipov-data}
  X_{n,i}(t) =
  \begin{cases}
  \varepsilon_{X,i}(t) \, , &i<k^* = \lfloor s^*n \rfloor +1 \\
  \varepsilon_{X,i}(t)(1+d_1+d_2(1+\sin(2\pi t))) \, , &i\geq k^*
  \end{cases}
\end{align}
where $\varepsilon_{X,1},\dots,\varepsilon_{X,n}$ are independent standard
Brownian motions. A sample size of $n=100$ is considered and the true change
point is defined by $k^* = 51$. The empirical rejection probabilities of the three tests are displayed in
Table~\ref{tab:sharipov}. The level ($d_1=d_2=0$) is approximated very well by all procedures under consideration.
Moreover,
the test   \eqref{eq:classical-test-cp} proposed in this paper is at least competitive in all
cases under consideration. In the case $d_1 = 0.4, d_2 = 0$ the procedures of
\cite{sharipov2019} perform slightly better but whenever $d_2 > 0$, the new
procedure  shows the best performance.

\begin{table}[t]
  \begin{center}
    \setlength{\tabcolsep}{5pt}
        { \scriptsize
    \begin{tabular}{c|ccc||c|ccc}
      $d_1, \, d_2$          & 1\%  & 5\%  & 10\% & $d_1, \, d_2$           &  1\% & 5\%          & 10\%  \\
      \hline \hline
      \multirow{2}{*}{0, 0}   & 1.3  & 5.0  & 9.9  & \multirow{2}{*}{0.4, 0}   & 19.3 & 44.5 & 61.0        \\
      & (0.4, 0.6)   & (4.4, 4.7)   & (10.0, 10.8) & & (16.1, 19.8) & (46.8, 50.1) & (63.2, 65.4) \\[4pt]
      \hline
      \multirow{2}{*}{0.8, 0} & 60.0 & 88.4 & 95.2 & \multirow{2}{*}{0, 0.4}   & 22.4 & 48.9 & 65.3         \\
      & (56.0, 58.8) & (88.0, 88.4) & (96.0, 95.5) & &  (9.8, 12.5) & (33.0, 38.6)  & (49.4, 55.4)  \\[4pt]
      \hline
      \multirow{2}{*}{0, 0.8} & 69.8 & 93.6 & 98.0 & \multirow{2}{*}{0.4, 0.4} &  63.4 & 89.3 & 95.8         \\
      & (45.8, 50.3) & (82.9, 85.8) & (93.8, 94.6) & & (44.2, 49.1) & (81.1, 82.2)  & (91.6, 92.1)  \\[4pt]
    \end{tabular}
    }
    \caption{\it Empirical rejection probabilities of the bootstrap   test \eqref{eq:classical-test-cp} for the classical hypotheses \eqref{class} of a structural break in the covariance operator. 
      The numbers in the brackets display the empirical rejection probabilities of the test proposed in \cite{sharipov2019}
      based on the supremum type   integral type CUSUM
      statistic (for $p=3$). 
      } \label{tab:sharipov}
  \end{center}
\end{table}

 Next we provide a comparison with the procedure proposed by \cite{stoehr2019}. Following these authors, we simulate  fAR(1) data where the errors (similar as in \eqref{eq:errors})
 are defined by
\begin{align*}
  e_j = \sum_{i=1}^{55} N_{i,j} \tilde \nu_i \, ,
  \qquad j=1,\dots , n \, ,
\end{align*}
 $\tilde \nu_1,\dots,\tilde\nu_{55}$ denote the Fourier basis and the random coefficients
\(N_{1,j},N_{2,j}, \dots, N_{55,j} \) are independent normally distributed with
expectation zero and variance $\var(N_{i,j})=\sigma_i^2$
($i=1,\ldots , 55$; $j=1,\ldots , n$). The fAR(1) data are then defined by
\begin{align} \label{eq:fAR1}
  X_{n,j} = \Psi(X_{n,j-1}) + e_j \, , \qquad j=1,\ldots,n \, ,
\end{align}
 where the linear operator $\Psi$ is represented by a
$55 \times 55$ matrix that is applied to the vector of the coefficients in the
basis representation. Here the matrix with $0.4$ on the diagonal and $0.1$ on
the superdiagonal and subdiagonal is chosen, such that  the generated fAR(1) time series is stationary. For the
  alternative  a change is inserted in the first $m$ leading
eigendirections for $m = 2,6,25$ by adding an additional normally distributed
noise term with variance $\sigma_\epsilon^2 / m$ for the observations $X_{n,j}$
for $j> \lfloor 0.5n \rfloor$. The following three settings are considered:
\begin{align*}
& \text{Setting 1: }~ \sigma_i = 1 ~\text{ for }~ i = 1,\dots,8
~\text{ and }~ \sigma_i = 0 ~\text{ for }~ i = 9,\dots, 55, ~~\sigma_\epsilon = 1.5 \\
& \text{Setting 2: }~ \sigma_i = 3^{-i} ~\text{ for }~ i = 1,\dots,55, ~~\sigma_\epsilon = 0.3 \\
& \text{Setting 3: }~ \sigma_i = i^{-1} ~\text{ for }~ i = 1,\dots,55, ~~\sigma_\epsilon = 1 \, .
\end{align*}

The empirical rejection probabilities of the test \eqref{eq:classical-test-cp} with block length $l=6$ and
the test based on dimension reduction developed in \cite{stoehr2019}  are displayed in   Table~\ref{tab:stoehr_power}.
We observe that in all cases under consideration the procedure proposed here yields an improvement with respect to the power.
 Note that
\cite{stoehr2019} also considered test procedures based on fully functional and
weighted functional statistics. As these methods considerably overestimate the test
level (see Figure~2 in the same reference),    these procedures are not included in the comparison.


\begin{table}[t]
  \centering  \setlength{\tabcolsep}{5.5pt}
      { \scriptsize
  \begin{tabular}{c||c|c|c}
    $m$ & Setting 1 & Setting 2 & Setting 3  \\
    \hline \hline
    0  & 4.7 (3.1)   & 8.1 (5.0)   & 3.9 (4.6) \\
    2  & 37.2 (22.8) & 92.5 (50.5) & 86.2 (30.4) \\
    6  & 81.1 (20.4) & 99.9 (98.9) & 99.9 (94.8) \\
    25 & 100 (29.0)  & 100 (92.3)  & 100 (97.3) \\
  \end{tabular}
  }
  \caption{\it Empirical rejection probabilities   (at level $5\%$) of the bootstrap test \eqref{eq:classical-test-cp} and the dimension reduction
    approach proposed in \cite{stoehr2019}
    (numbers in the brackets).  } \label{tab:stoehr_power}
\end{table}

\paragraph{Relevant hypotheses:}

We conclude   this section   investigating the finite sample properties of the test defined by \eqref{eq:relevant-test-cp} for the  hypotheses \eqref{eq:relevant-hypothesis}
of a relevant change in the covariance operator. For this purpose we consider similar scenarios  as in Section~\ref{sec:simulations-relevant}. In all cases, the location of the change  is set to $s^* = 0.5$ and the observations after the change point are multiplied by a constant $a$ such that \eqref{h00} holds.
For the estimation of the extremal sets, the  parameter in \eqref{eq:estimatedSets-cp} is set as $c_n = 0.1 \log(n)$.

In Table~\ref{tab:relevant-lvl-cp} empirical rejection probabilities are displayed for different processes at the boundary of  the null hypothesis i.e. the observations after the change point are multiplied by $a = 2$  and the threshold $\Delta$ is defined in each case such that  $\| C_1 - C_2\|_\infty=\Delta$. For fIID and non-Gaussian data
the block length in \eqref{eq:bootProcess-cp} is set to $l=1$ and the threshold  is given by $\Delta = a^2-1$. The fMA(1) and fMA(2) data are  defined by \eqref{eq:fMA2} with $\kappa_1 = 0.7, \kappa_2 = 0$ and $\kappa_1 = 0.5, \kappa_2 = 0.3$, respectively. The threshold parameter is set to $\Delta = (a^2-1)(1+\kappa_1^2+\kappa_2^2)$ in both cases and the block length in \eqref{eq:bootProcess-cp} is set to $l=2$ and $l=3$, respectively.

We observe that the nominal   level is reasonably well approximated in most cases under consideration especially for the sample size $n=200$. Only in the non-Gaussian case, the nominal level is underestimated for the sample size $N=100$, but the approximation improves considerably for the sample size $N=200$.

In Table~\ref{tab:relevant-power-cp}, we show the empirical rejection
probabilities  of the test \eqref{eq:relevant-test-cp} and also for the test
developed in \cite{detkokvol2020} for scenarios in the interior of the null
hypothesis of no relevant change point as well as under the alternative. We
consider independent identically distributed Gaussian (fIID) and fMA(2) data
and multiply the observations after the change point by different values
$a = 1.8,1.9,2,2.2,2.4,2.6$. In the fIID case the threshold parameter is given by   $\Delta = 3$ and in the fMA(2) case it
is   $\Delta =  3\cdot(1+\kappa_1^2+\kappa_2^2)$ (where still
$\kappa_1 = 0.5, \kappa_2 = 0.3$). Consequently, the case $a=2$  always corresponds to
the boundary of the null hypothesis, and the cases $a<2$ and $a >2$ represent the interior of the    null hypothesis and
alternative.
Since the procedure developed by \cite{detkokvol2020} is based on a different metric, the threshold parameter $\Delta$ in the relevant hypotheses \eqref{eq:relevant-intro} is set to
\begin{align*}
\Delta = \int_{[0,1]} \int_{[0,1]} \{(1 - 2^2)C_1(s,t)\}^2 ds dt
\end{align*}
for this test procedure. Consequently the boundary of the null hypothesis of no relevant change in the covariance operators (w.r.t. the corresponding metric) is also obtained for the factor $a=2$ for both data models.

We mention again that the nominal level at the boundary of the hypotheses is
reasonably well approximated by the test \eqref{eq:relevant-test-cp} while the
test procedure developed in \cite{detkokvol2020} is more conservative. In the
interior of the null hypothesis $(a < 2)$ the rejection probabilities of both
tests are strictly smaller than the nominal level. This property is desirable
as it means that the probability of a type I error is small in situations with
a large deviation from the alternative. On the other hand,  under  the alternative
 the new test \eqref{eq:relevant-test-cp} has
substantially more power than the test developed in \cite{detkokvol2020}.

\begin{table}[t]
  \centering  \setlength{\tabcolsep}{5.5pt}
      { \scriptsize
  \begin{tabular}{c|ccc||ccc||ccc||ccc}
    & \multicolumn{3}{|c||} {fIID} & \multicolumn{3}{|c||} {non-Gaussian} & \multicolumn{3}{|c||} {fMA(1)} & \multicolumn{3}{|c} {fMA(2)} \\
    \hline
    \(n\)     & 1\% & 5\% & 10\% & 1\% & 5\% & 10\% & 1\% & 5\% & 10\% & 1\% & 5\% & 10\%  \\
    \hline \hline
    100      & 1.1 & 3.8 & 9.5  & 0   & 0.8 & 5.3  & 0.8 & 4.9 & 13.7 & 1.3 & 6.0 & 11.7  \\
    200      & 0.7 & 4.6 & 10.1 & 0.3 & 3.1 & 8.4  & 1.3 & 4.9 & 9.8  & 0.7 & 4.9 & 10.5
  \end{tabular}
  }
  \caption{\it \label{tab:relevant-lvl-cp}
    Simulated level of the test \eqref{eq:relevant-test-cp} for the hypotheses
    defined by \eqref{eq:relevant-hypothesis} at the boundary of the hypotheses, that is $\|C_1-C_2\|_\infty=\Delta$.
     }
\end{table}

\begin{table}[h]
  \centering
      { \scriptsize
  \begin{tabular}{c||c|c||c|c }
    & \multicolumn{2}{c||} {fIID} & \multicolumn{2}{c}{fMA(2)}  \\
    \cline{1-5}
    $a$ & $n=100$ & $n=200$ & $n=100$ & $n=200$ \\
    \hline \hline
    1.8 & 0.3 (0.4)   & 0 (0)       & 1.3 (0.1)   & 1.0 (0.1)   \\
    1.9 & 1.8 (0.9)   & 0.1 (0.5)   & 3.4 (0.4)   & 1.4 (0.4)   \\
    2.0 & 3.8 (2.3)   & 4.6 (3.2)   & 6.0 (1.0)   & 4.9 (1.4)   \\
    2.2 & 21.5 (9.8)  & 33.6 (25.4) & 19.4 (6.2)  & 27.2 (11.1) \\
    2.4 & 47.0 (23.3) & 74.9 (51.2) & 40.0 (15.6) & 65.9 (31.2) \\
    2.6 & 73.0 (37.9) & 96.0 (70.6) & 63.3 (26.5) & 88.0 (49.7)
  \end{tabular}
  }
  \caption{\it \label{tab:relevant-power-cp}
    Simulated rejection probabilities  of the test \eqref{eq:relevant-test-cp}
    for the hypotheses   \eqref{eq:relevant-hypothesis} of a relevant change in the covariance operator
    considering fIID
     and fMA(2) data (level $5\%$).  The cases $a<2, a=2$ and $a>2$ correspond to the interior, boundary of the null hypothesis and to the alternative. The numbers in   brackets represent the empirical rejection probabilities of the procedure developed in \cite{detkokvol2020}.
  }
\end{table}

\subsection{Data Example} \label{sec52}

Similar as \cite{fremdt2013} and \cite{Paparoditis2016} we consider egg-laying
curves of medflies (Mediterranean fruit flies, Ceratitis capitata). The
original data consists of the number of eggs which were laid on each day during
the lifetime of $1000$ female medflies and a detailed description of the
experiment can be found in \cite{carey1998}. Only medflies which lived at least
$34$ days are considered and split into two samples, the medflies which lived
at most $43$ days and those which lived at least $44$ days. A Fourier basis
consisting of $49$ basis functions is used to transform the discrete
observations to functional data $(X_j\colon j = 1,\dots,256)$ and
$(Y_j\colon j = 1,\dots,278)$. The expressions $X_i(t)$ and $Y_j(t)$ denote the
number of eggs which were laid on day $\lfloor 30 t \rfloor$ by the $i$th
short-lived and the $j$th long-lived medfly relative to the total number of
eggs laid in the whole lifetime of the $i$th short-lived and the $j$th
long-lived medfly, respectively ($t\in [0,1], i = 1,\dots,256,j = 1,\dots,278$).
First, the test \eqref{eq:classical-boottest} is used to study the classical
hypotheses in \eqref{class}. The window parameters in \eqref{eq:bootProcess}
are set to $l_1=l_2=1$ since the egg-laying curves corresponding to the
different medflies can be regarded as independent. For the calculation of
critical values, $200$ bootstrap samples are generated. The classical null hypothesis of
equal covariance operators is then rejected at level $5\%$ and can not be
rejected at level $1\%$. The outcome when using the procedure developed in
\cite{fremdt2013} depends on the choice of the number of considered functional
principal components $p$ and the procedure developed in \cite{Paparoditis2016}
yields a $p$-value of $0.3\%$ (see Table~3 in \cite{Paparoditis2016}).
In Table~\ref{tab:data-example} the empirical rejection probabilities of the
test \eqref{eq:relevant-boottest} are displayed for the relevant hypotheses in
\eqref{eq:relevant-hypothesis} for different choices of the threshold parameter
$\Delta$.  It can be seen that even for $\Delta = 0.0003$
i.e. when a maximal deviation of only $0.0003$ is tolerated, the null
hypothesis of no relevant difference between the covariance operators can not
be rejected at all considered test levels. For $\Delta = 0.0002$ the null can
be rejected at level $10\%$ and for $\Delta = 0.0001$ also at level $5\%$.
Although the classical null hypothesis of equal covariance operators is
rejected at level $5\%$, these results may raise the question if the detected
difference is really of practical relevance.

\begin{table}[h]
  \vspace{.5cm}
  \begin{center}
        { \scriptsize
    \begin{tabular}{c@{\qquad}lll}
      \hline
      $\Delta$ & 1\%          & 5\%         & 10\%  \\
      \hline
      0.0001 & FALSE  & TRUE  & TRUE  \\
      0.0002 & FALSE & FALSE  & TRUE  \\
      0.0003 & FALSE & FALSE  & FALSE  \\
      \hline
    \end{tabular}
    }
  \end{center}
  \caption{\it Summary of the outcome of the test \eqref{eq:relevant-boottest} for the relevant hypotheses \eqref{eq:relevant-hypothesis} for different values of $\Delta$ for the relative egg-laying curves of medflies. The label TRUE means that the null hypothesis is rejected and the label FALSE means that the null hypothesis is not rejected.}
  \label{tab:data-example}
\end{table}




\section{Appendix: Proofs of main results} \label{sec7}
\def\theequation{6.\arabic{equation}}
\setcounter{equation}{0}

\subsection{Proof of Theorem~\ref{thrm:CLT}} \label{sec61}
We apply the central limit theorem as formulated in Theorem~2.1 in
\cite{dette2018} to the sequence of $C(T^2)$-valued random variables
$((Z_{j} - \mu)^{\check \otimes 2})_{j\in\N} = (\eta_{j}^{\check \otimes 2})_{j\in\N}$.

It can be easily seen that conditions (A1), (A2) and (A4) in this reference
are satisfied.
In order to see that the remaining condition (A3) also holds, we use the triangle
inequality and    Assumption~\ref{assumption} of the present work
to obtain, for any $j\in\N$ and $s,t,s^\prime,t^\prime\in T$,
\begin{align*}
  |\eta_{j}(s)\eta_{j}(t) - \eta_{j}(s^\prime)\eta_{j}(t^\prime)|
  &\leq |\eta_{j}(s)(\eta_{j}(t) - \eta_{j}(t^\prime))|
  + |\eta_{j}(t^\prime)(\eta_{j}(s) - \eta_{j}(s^\prime))| \\
  &\leq \|\eta_{j}\|_\infty \, \big(|\eta_{j}(t) - \eta_{j}(t^\prime)|
  + |\eta_{j}(s) - \eta_{j}(s^\prime)| \big) \\
  &\leq \|\eta_{j}\|_\infty \, M \, \big( \rho(t, t^\prime)
  + \rho(s, s^\prime) \big) \\
  &\lesssim \|\eta_{j}\|_\infty \, M \,
  \rho_{\max} \big((t,s), (t^\prime, s^\prime) \big)
\end{align*}
where $\E \big[ (\|\eta_{j}\|_\infty \, M )^J \big] \leq \tilde K < \infty$ by (A3).
Now observe that
\begin{align*}
  \frac{1}{\sqrt{n}} \sum_{j=1}^n (Z_{j} - \bar{Z}_{n})^{\check \otimes 2}
  = \frac{1}{\sqrt{n}} \sum_{j=1}^n \eta_{j}^{\check \otimes 2}
  -\frac{1}{\sqrt{n}}
  \bigg(\frac{1}{\sqrt{n}} \sum_{j=1}^n \eta_{j} \bigg)^{\check \otimes 2}
  = \frac{1}{\sqrt{n}} \sum_{j=1}^n \eta_{j}^{\check \otimes 2} + o_\P(1)
\end{align*}
which yields the claim since Theorem~2.1 in \cite{dette2018} can be applied to
the sequence $( \eta_{j}^{\check \otimes 2} )_{j\in\N}$ as shown above.

\subsection{Proof of Proposition \ref{prop1}} \label{sec62}

As the samples are independent, it directly
follows from Theorem~\ref{thrm:CLT} that
\begin{align*}
  \sqrt{m+n}& \bigg (
  \frac{1}{m}
  \sum_{j=1}^m (\tilde X_{m,j}^{\check \otimes 2} - C_1 ),
  \, \frac{1}{n}
  \sum_{j=1}^n (\tilde Y_{n,j}^{\check \otimes 2} - C_2 ) \bigg ) \\
  &= \sqrt{m+n} \bigg (
  \frac{1}{m}
  \sum_{j=1}^m (\eta_{1,j}^{\check \otimes 2} - C_1 ),
  \, \frac{1}{n}
  \sum_{j=1}^n (\eta_{2,j}^{\check \otimes 2} - C_2 ) \bigg ) + o_\P(1)
  \rightsquigarrow  \bigg ( \frac{1}{\sqrt{\lambda}} ~ Z_1,
  \frac{1}{\sqrt{1-\lambda}} ~ Z_2 \bigg )
\end{align*}
in $C([0,1]^2)^2$ as $m,n\to\infty$, where $Z_1$ and $Z_2$ are independent,
centred Gaussian processes defined by their long-run covariance operators
\eqref{c1} and \eqref{c2}.
By the continuous mapping theorem it follows that
\begin{align} \label{eq:Z-limit}
  Z_{m,n} = \sqrt{m+n} \, \bigg (
  \frac{1}{m}
  \sum_{j=1}^m \tilde X_{m,j}^{\check \otimes 2} -
  \, \frac{1}{n}
  \sum_{j=1}^n \tilde Y_{n,j}^{\check \otimes 2}  - (C_1 - C_2) \bigg )
  \rightsquigarrow  Z
\end{align}
in $C([0,1]^2)$ as $m,n\to\infty$,  where $Z$ is again a centred Gaussian
process with covariance operator  \eqref{eq:Z-cov}.

If $d_\infty = 0$, the convergence in \eqref{eq:Z-limit} together with the
continuous mapping yield  \eqref{eq:sup_T}.
If $d_\infty > 0$, the asymptotic distribution of $\hat d_\infty$ can be
deduced from Theorem~B.1 in the online supplement of \cite{dette2018} or  alternatively
from the results in \cite{carcamo2019}.

\subsection{Proof of Theorem~\ref{thrm:bootstrap_class}  and  \ref{thrm:bootstrap_test}}

{\bf Proof of Theorem ~\ref{thrm:bootstrap_class}.}
Using similar  arguments as in the proof of Theorem~\ref{thrm:CLT}, it  follows that
the process $ \hat B^{(r)}_{m,n}$ in \eqref{eq:bootProcess} admits the stochastic expansion
\begin{align*}
  \hat B^{(r)}_{m,n} = \sqrt{n+m} \bigg\{
  &\frac{1}{m} \sum_{k=1}^{m-l_1+1} \frac{1}{\sqrt{l_1}}\bigg(
  \sum_{j=k}^{k+l_1-1} \eta_{1,j}^{\check \otimes 2}
  -\frac{l_1}{m}\sum_{i=1}^m \eta_{1,j}^{\check \otimes 2}  \bigg)
  \xi_k^{(r)} \\
  - &\frac{1}{n} \sum_{k=1}^{n-l_2+1} \frac{1}{\sqrt{l_2}}\bigg(
  \sum_{j=k}^{k+l_2-1} \eta_{2,j}^{\check \otimes 2}
  -\frac{l_2}{n}\sum_{i=1}^n \eta_{2,j}^{\check \otimes 2} \bigg)
  \zeta_k^{(r)} \bigg\}  + o_\P(1) \, ,
\end{align*}
and the sequences $(\eta_{1,j}^{\check \otimes 2})_{j\in\N}$ and
$(\eta_{2,j}^{\check \otimes 2})_{j\in\N}$ satisfy Assumption~2.1 in
\cite{dette2018}. Thus, similar arguments as in the proof of Theorem~3.3 in the same reference yield
\begin{align} \label{eq:vec-Z}
\big( Z_{m,n} ,
\hat B_{m,n}^{(1)},\dots,\hat B_{m,n}^{(R)}\big)
\rightsquigarrow
(Z, Z^{(1)},\dots,Z^{(R)})
\end{align}
in $C([0,1]^2)^{R+1}$ as $m,n \to \infty$ where the process $Z_{m,n}$ is defined in \eqref{eq:Z-limit} and the random functions $Z^{(1)},\dots,Z^{(R)}$ are independent copies of $Z$ which is also defined in  \eqref{eq:Z-limit}.

If $d_\infty = 0$, the continuous mapping theorem implies
\begin{align} \label{eq:vec-T}
\big( \sqrt{m+n} \, \hat{d}_{\infty} ,~
T_{m,n}^{(1)},\dots,T_{m,n}^{(R)}\big)
\stackrel{\mathcal{D}}{\longrightarrow} (T,~ T^{(1)},\dots,T^{(R)})
\end{align}
in $\R^{R+1}$ as $m,n \to \infty$ where the statistic $\hat{d}_\infty$ is
defined by \eqref{eq:d-hat}, the bootstrap statistics $T_{m,n}^{(1)},\dots,T_{m,n}^{(R)}$
are defined by \eqref{eq:bootStat-classic} and the random variables
$T^{(1)},\dots,T^{(R)}$ are independent copies of $T$ which is defined by
\eqref{eq:sup_T}. Now, Lemma~4.2 in \cite{buecher2019} directly implies \eqref{eq:boot_alpha_classic}, that is
\begin{align*}
\lim_{m,n,R\to\infty} \mathbb{P}\bigg ( \hat{d}_{\infty}
>  \frac{T_{m,n}^{\{\lfloor R(1-\alpha)\rfloor\}}}{\sqrt{m+n}} \bigg)
= \alpha \, .
\end{align*}
For the application of this result, it is required that the distribution of the
random variable $T$ has a continuous distribution function,  which follows from \cite{Gaenssler2007}.
In order to show the consistency of test \eqref{eq:classical-boottest} in the case $d_\infty >0$, write
\begin{align*}
  \mathbb{P}\bigg ( \hat{d}_{\infty}
  >  \frac{T_{m,n}^{\{\lfloor R(1-\alpha)\rfloor\}}}{\sqrt{m+n}} \bigg)
  &= \mathbb{P}\big(
  \sqrt{m+n} \, (\hat{d}_{\infty} - d_{\infty}) + \sqrt{m+n} \, d_{\infty}
  > T_{m,n}^{\{\lfloor R(1-\alpha)\rfloor\}} \big)
\end{align*}
and note that, given \eqref{eq:vec-T} and \eqref{eq:T-limit}, the assertion in \eqref{eq:consistency} follows by simple arguments.

\medskip

{\bf Proof of Theorem \ref{thrm:bootstrap_test}.}
  First
note that the same arguments as in the proof of Theorem~3.6 in \cite{dette2018}
show that the estimators of the extremal sets defined by
\eqref{eq:estimatedSets} are consistent that is
\begin{align*}
d_H( \hat{\mathcal{E}}_{m,n}^\pm ,  \mathcal{E}^\pm)
\xrightarrow[m,n\to\infty]{\P}  0 \, ,
\end{align*}
where $d_H$ denotes the Hausdorff distance.
Thus, given the convergence in \eqref{eq:vec-Z}, the  arguments  in the proof of Theorem~3.7 in the same reference yield
\begin{align} \label{eq:vec-convergence}
  \big( \sqrt{n+m} ~ (\hat{d}_\infty - d_\infty) ,~
  K_{m,n}^{(1)},\dots,K_{m,n}^{(R)}\big)
  \stackrel{\mathcal{D}}{\longrightarrow}
  (T(\mathcal{E}),~ T^{(1)}(\mathcal{E}),\dots,T^{(R)}(\mathcal{E}))
\end{align}
in $\R^{R+1}$ as $m,n \to \infty$ where the statistic $\hat{d}_\infty$ is
defined by \eqref{eq:d-hat}, the bootstrap statistics
$K_{m,n}^{(1)},\dots,K_{m,n}^{(R)}$ are defined by \eqref{eq:bootStat} and the
random variables $T^{(1)}(\cal E),\dots,T^{(R)}(\cal E)$ are independent copies
of $T(\mathcal{E})$ which is defined by \eqref{eq:T-limit}. Note that this
convergence holds true under the null and the alternative hypothesis.

If $\Delta = d_\infty$, Lemma~4.2 in \cite{buecher2019} directly implies \eqref{eq:boot_alpha} and again the results in \cite{Gaenssler2007} ensure that the limit $T(\mathcal{E})$ has a continuous distribution function.

If $\Delta \neq d_\infty$, write
\begin{align*}
 \mathbb{P}\bigg ( \hat{d}_{\infty}
> \Delta + \frac{K_{m,n}^{\{\lfloor R(1-\alpha)\rfloor\}}}{\sqrt{n+m}} \bigg)
= \mathbb{P}\big(
\sqrt{m+n} \, (\hat{d}_{\infty} - d_{\infty}) + \sqrt{m+n} \, (d_{\infty} - \Delta)
> K_{m,n}^{\{\lfloor R(1-\alpha)\rfloor\}} \big) \, .
\end{align*}
Then it follows from \eqref{eq:vec-convergence} and simple arguments that, for any
$R\in\N$,
\begin{align*}
\lim_{m,n\to\infty} \mathbb{P}\bigg ( \hat{d}_{\infty}
> \Delta + \frac{K_{m,n}^{\{\lfloor R(1-\alpha)\rfloor\}}}{\sqrt{n+m}} \bigg)
= 0 \quad \text{and} \quad
\liminf_{m,n\to\infty} \mathbb{P}\bigg ( \hat{d}_{\infty}
> \Delta + \frac{K_{m,n}^{\{\lfloor R(1-\alpha)\rfloor\}}}{\sqrt{n+m}} \bigg)
= 1
\end{align*}
if $\Delta > d_\infty$ and $\Delta < d_\infty$, respectively. This proves the
remaining assertions of Theorem~\ref{thrm:bootstrap_test}.

\subsection{Proof of Proposition \ref{prop2}}  \label{sec64}

Let $C_{n,j}$ denote the covariance operator of $X_{n,j}$ defined by $C_{n,j}(s,t) = \cov (X_{n,j}(s),X_{n,j}(t))$
and consider the sequential process
\begin{align*}
  \hat{\mathbb{V}}_n (s) &  = \frac{1}{\sqrt{n}} \sum^{\lfloor s  n \rfloor}_{j=1}
  (\tilde X_{n,j}^{\check \otimes 2}  - C_{n,j})
  + \sqrt{n} \Big (s \, - \frac {\lfloor s  n \rfloor}{n} \Big )
  \big(\tilde X_{n, \lfloor s  n \rfloor+1}^{\check \otimes 2}
  - C_{n,\lfloor s  n \rfloor+1} \big) \,  \\
&  =  \frac{1}{\sqrt{n}} \sum^{\lfloor s  n \rfloor}_{j=1}
  (\tilde \eta_{n,j}^{\check \otimes 2}  - C_{n,j})
  + \sqrt{n} \Big (s \, - \frac {\lfloor s  n \rfloor}{n} \Big )
  \big(\tilde \eta_{n, \lfloor s n \rfloor+1}^{\check \otimes 2}
  - C_{n,\lfloor s  n \rfloor+1} \big) + o_\P(1)
\end{align*}
which is an element of $C([0,1], C([0,1]^2))$. Note that $\{\hat{\mathbb{V}}_n(s)\}_{s\in[0,1]}$ can equivalently be regarded as an
element of $C([0,1]^3)$ and we have the representation
\begin{align} \label{eq:Vn-representation}
  \hat{\mathbb{V}}_n
  = \tilde{\mathbb{V}}_{1,n} + \tilde{\mathbb{V}}_{2,n} \, ,
\end{align}
where the processes $ \tilde{\mathbb{V}}_{1,n},  \tilde{\mathbb{V}}_{2,n} \in C([0,1]^3)$ are defined by
\begin{align*}
   \tilde{\mathbb{V}}_{1,n}(s,t,u) &= \, \hat{\mathbb{V}}_{1,n}(s,t,u) \mathds{1}\{\lfloor sn \rfloor < \lfloor s^* n \rfloor\}
   + \hat{\mathbb{V}}_{1,n}(\lfloor s^* n \rfloor/n,t,u) \mathds{1}\{\lfloor sn \rfloor \geq \lfloor s^* n \rfloor\} \\
   \tilde{\mathbb{V}}_{2,n}(s,t,u) &= \, (\hat{\mathbb{V}}_{2,n}(s,t,u) - \hat{\mathbb{V}}_{2,n}(\lfloor s^* n \rfloor/n,t,u))
   \mathds{1}\{\lfloor sn \rfloor \geq \lfloor s^* n \rfloor\}
\end{align*}
($s,t,u\in[0,1]$) and
\begin{align*}
\hat{\mathbb{V}}_{l,n}(s) =  \frac{1}{\sqrt{n}} \sum^{\lfloor s   n \rfloor}_{j=1}
(\eta_{l,j}^{\check \otimes 2}  - C_{l})
+ \sqrt{n} \Big (s \, - \frac {\lfloor s  n \rfloor}{n} \Big )
\big(\eta_{l, \lfloor s  n \rfloor+1}^{\check \otimes 2}
- C_{l} \big) \quad \quad (l = 1,2) \, .
\end{align*}

Recall the definition of the array ($\tilde{\eta}_{n,j} \colon n\in\N, j = 1,\dots, n$) in \eqref{eq:tilde-eta}.
By Theorem~2.2 in \cite{dette2018} it follows that
\begin{align*}
  \hat{\mathbb{V}}_{l,n}   \rightsquigarrow \mathbb{V}_l \quad \quad (l = 1,2)
\end{align*}
in $C([0,1]^3)$, where  $\mathbb{V}_l$ is a centred Gaussian measure on
$C([0,1]^3)$ characterized by the covariance operator
\begin{align*}
  \cov \big (\mathbb{V}_l(s,t,u), \mathbb{V}_l(s^\prime,t^\prime,u^\prime) \big)
  &= (s \wedge s^\prime) \, \mathbb{C}_l((t,u),(t^\prime, u^\prime)), \qquad l = 1,2
\end{align*}
and the long-run covariance operator $\mathbb C_l$ is defined in \eqref{eq:lrv-cp}. From the continuous mapping theorem we obtain
\begin{align} \label{eq:WIP-l}
\tilde{\mathbb{V}}_{l,n}   \rightsquigarrow \tilde{\mathbb{V}}_l \quad \quad (l = 1,2)
\end{align}
in $C([0,1]^3)$, where  $\tilde{\mathbb{V}}_1, \tilde{\mathbb{V}}_2$ are
centred Gaussian measures on $C([0,1]^3)$ characterized by
\begin{align*}
\tilde{\mathbb{V}}_1 (s,t,u) = \mathbb{V}_{1}(s \wedge s^*,t,u) \, , \quad
\tilde{\mathbb{V}}_2 (s,t,u) = (\mathbb{V}_{2}(s,t,u) - \mathbb{V}_{2}(s^*,t,u))
\mathds{1}\{ s \geq s^* \}
\end{align*}
with covariance operators
\begin{align*}
\cov \big (\tilde{\mathbb{V}}_1(s,t,u), \tilde{\mathbb{V}}_1(s^\prime,t^\prime,u^\prime) \big)
&= (s \wedge s^\prime \wedge s^*) \, \mathbb{C}_1((t,u),(t^\prime, u^\prime)) \\
\cov \big (\tilde{\mathbb{V}}_2(s,t,u), \tilde{\mathbb{V}}_2(s^\prime,t^\prime,u^\prime) \big)
&= (s\wedge s^\prime - s^*)_+ \, \mathbb{C}_2((t,u),(t^\prime, u^\prime)) \, .
\end{align*}
In the following we will show the weak convergence
\begin{align} \label{eq:WIP}
  \hat{\mathbb{V}}_n \rightsquigarrow \mathbb{V}
\end{align}
in $C([0,1]^3)$ as $n\to\infty$, where $\mathbb{V}\in C([0,1]^3)$ is a centred Gaussian random variable characterized by its covariance operator
\begin{align*}
  \cov(\mathbb{V}(s,t,u), \mathbb{V}(s^\prime,t^\prime,u^\prime)) = (s\wedge s^\prime \wedge s^*) \, \mathbb{C}_1((t,u), (t^\prime, u^\prime)) + (s\wedge s^\prime - s^*)_+ \, \mathbb{C}_2((t,u), (t^\prime, u^\prime))
\end{align*}
and the long-run covariance operators $\mathbb{C}_1, \mathbb{C}_2$
are defined by \eqref{eq:lrv-cp}. The convergence in \eqref{eq:WIP-l} implies
that the processes $\tilde{\mathbb{V}}_{1,n}, \tilde{\mathbb{V}}_{2,n}$ are asymptotically tight and the representation in \eqref{eq:Vn-representation} yields that $\hat{\mathbb{V}}_{n}$ is asymptotically tight as well \citep[see Section~1.5 in][]{wellner1996}.
In order to prove the convergence in \eqref{eq:WIP} it consequently remains to
show the convergence of the finite dimensional distributions. For this, we utilize the Cr\'amer-Wold device and show that
\begin{align*}
  \tilde Z_n =\sum_{j=1}^q c_j \hat{\mathbb{V}}_n (s_j,t_j,u_j)
 & =  \sum_{j=1}^q c_j \big\{ \tilde{\mathbb{V}}_{1,n}(s_j,t_j,u_j) +
  \tilde{\mathbb{V}}_{2,n}(s_j,t_j,u_j) \big\}  \\
 &  \stackrel{\cal D}{\longrightarrow} \tilde Z = \sum_{j=1}^q c_j  \mathbb{V}(s_j ,t_j,u_j)
\end{align*}
for any $(s_1,t_1,u_1),\dots,(s_q,t_q,u_q) \in [0,1]^3$, $c_1,\dots,c_q \in \R$ and $q\in\N$. Asymptotic normality of $\tilde{Z}_n$ can be proved by the same arguments as in the proof of Theorem~2.1 in \cite{dette2018} and it remains to show that the variance of the random variable $\tilde Z_n$ converges to the variance of $\tilde Z$. Using (3.17) in \cite{dehling2002} and assumptions (A2) and (A4) we obtain for any $(s,t,u),(s^\prime, t^\prime, u^\prime) \in [0,1]^3$
\begin{align} \label{eq:cov0}
\begin{split}
  &\cov (\tilde{\mathbb{V}}_{1,n} (s,t,u), \tilde{\mathbb{V}}_{2,n} (s^\prime,t^\prime,u^\prime) ) \\
  &= \frac{1}{n} \sum^{\lfloor (s \wedge s^*) n \rfloor}_{j=1} \sum^{\lfloor s^\prime n \rfloor}_{i=\lfloor s^* n \rfloor + 1}
  \cov (\tilde{\eta}_{n,j}^{\check \otimes 2}(t,u), \tilde{\eta}_{n,i}^{\check \otimes 2}(t^\prime,u^\prime)) + o(1) \\
  &\lesssim \frac{1}{n} \sum^{\lfloor (s \wedge s^*) n \rfloor}_{j=1} \sum^{\lfloor s^\prime n \rfloor}_{i=\lfloor s^* n \rfloor + 1} \|\tilde{\eta}_{n,j}^{\check \otimes 2}(t,u)\|_2 \, \|\tilde{\eta}_{n,i}^{\check \otimes 2}(t^\prime,u^\prime)\|_2 \, \varphi(i-j)^{1/2} + o(1) \\
  &\lesssim \frac{1}{n} \sum^{\lfloor (s \wedge s^*) n \rfloor}_{j=1} \sum^{\lfloor s^\prime n \rfloor}_{i=\lfloor s^* n \rfloor + 1}  \varphi(i-j)^{1/2} + o(1) \\
  &\lesssim \frac{1}{n} \sum^{\lfloor s^\prime n \rfloor - 1}_{i = 1} i \varphi(i)^{1/2} + o(1) \stackrel[n\to\infty]{}{\longrightarrow} 0 \, ,
\end{split}
\end{align}
where the symbol ``$\lesssim$'' means less or equal up to a constant independent of $n$, and $\|X\|_2 = \E [X^2]^{1/2}$ denotes the $L^2$-norm of a real valued random variable $X$ (also note that we implicitly assume $\sum_{i=j}^k a_i = 0$ if $k<j$).  Furthermore, assuming without loss of generality that $s \leq s^\prime$, we have
\begin{align*}
  &\cov (\tilde{\mathbb{V}}_{1,n} (s,t,u), \tilde{\mathbb{V}}_{1,n} (s^\prime,t^\prime,u^\prime) ) \\
  &= \frac{1}{n} \sum^{\lfloor (s\wedge s^*) n \rfloor}_{j=1} \sum^{\lfloor (s^\prime\wedge s^*) n \rfloor}_{i= 1}
  \cov (\eta_{1,j}^{\check \otimes 2}(t,u), \eta_{1,i}^{\check \otimes 2}(t^\prime,u^\prime)) + o(1) \\
  &= \frac{1}{n} \sum^{\lfloor (s\wedge s^*) n \rfloor}_{j=1} \bigg(
  \sum^{\lfloor (s \wedge s^*) n \rfloor}_{i= 1} + \sum^{\lfloor (s^\prime\wedge s^*) n \rfloor}_{i=\lfloor (s \wedge s^*) n \rfloor + 1} \bigg)
  \cov (\eta_{1,j}^{\check \otimes 2}(t,u), \eta_{1,i}^{\check \otimes 2}(t^\prime,u^\prime))  + o(1) \\
  &= \frac{1}{n} \sum^{\lfloor (s\wedge s^*) n \rfloor}_{j=1}
  \sum^{\lfloor (s \wedge s^*) n \rfloor}_{i= 1}
  \cov (\eta_{1,j}^{\check \otimes 2}(t,u), \eta_{1,i}^{\check \otimes 2}(t^\prime,u^\prime))  + o(1) \, ,
\end{align*}
where the  last equality follows by the same arguments as used in \eqref{eq:cov0}. For the remaining expression we use the dominated convergence theorem to obtain
\begin{align*}
  &\frac{1}{n} \sum^{\lfloor (s\wedge s^*) n \rfloor}_{j=1}
  \sum^{\lfloor (s \wedge s^*) n \rfloor}_{i= 1}
  \cov (\eta_{1,j}^{\check \otimes 2}(t,u), \eta_{1,i}^{\check \otimes 2}(t^\prime,u^\prime)) \\
  &= \sum^{\lfloor (s \wedge s^*) n \rfloor -1}_{i= -(\lfloor (s \wedge s^*) n \rfloor -1)}
  \frac{\lfloor (s \wedge s^*) n \rfloor - |i|}{n} \, \cov (\eta_{1,0}^{\check \otimes 2}(t,u), \eta_{1,i}^{\check \otimes 2}(t^\prime,u^\prime))
  \stackrel[n\to\infty]{}{\longrightarrow} (s \wedge s^*) \, \mathbb{C}_1((t,u), (t^\prime, u^\prime))
\end{align*}
which means that for any $(s,t,u),(s^\prime, t^\prime, u^\prime) \in [0,1]^3$
\begin{align*}
  \cov (\tilde{\mathbb{V}}_{1,n} (s,t,u), \tilde{\mathbb{V}}_{1,n} (s^\prime,t^\prime,u^\prime) )
  \stackrel[n\to\infty]{}{\longrightarrow} (s \wedge s^\prime \wedge s^*) \, \mathbb{C}_1((t,u), (t^\prime, u^\prime)) \, .
\end{align*}
By similar arguments we obtain
\begin{align*}
\cov (\tilde{\mathbb{V}}_{2,n} (s,t,u), \tilde{\mathbb{V}}_{2,n} (s^\prime,t^\prime,u^\prime) )
\stackrel[n\to\infty]{}{\longrightarrow} (s \wedge s^\prime -  s^*)_+ \, \mathbb{C}_2((t,u), (t^\prime, u^\prime))
\end{align*}
and therefore we have
\begin{align*}
  \var (\tilde Z_n ) &= \sum_{j=1}^q \sum_{j^\prime=1}^q c_j c_{j^\prime} \cov (\hat{\mathbb{V}}_n (s_j,t_j,u_j), \hat{\mathbb{V}}_n (s_{j^\prime},t_{j^\prime},u_{j^\prime}) ) \\
  &= \sum_{j=1}^q \sum_{j^\prime=1}^q c_j c_{j^\prime} \big\{ \cov (\tilde{\mathbb{V}}_{1,n} (s_j,t_j,u_j), \tilde{\mathbb{V}}_{1,n} (s_{j^\prime},t_{j^\prime},u_{j^\prime}) )  \\
  &\hspace{70pt}+ \cov (\tilde{\mathbb{V}}_{2,n} (s_j,t_j,u_j), \tilde{\mathbb{V}}_{2,n} (s_{j^\prime},t_{j^\prime},u_{j^\prime}) ) \big\}
  + o(1) \\
  &\stackrel[n\to\infty]{}{\longrightarrow}
  \sum_{j=1}^q \sum_{j^\prime=1}^q c_j c_{j^\prime}  \cov (\mathbb{V} (s_j,t_j,u_j), \mathbb{V} (s_{j^\prime},t_{j^\prime},u_{j^\prime}) )
  = \var(\tilde Z)
\end{align*}
which finally proves \eqref{eq:WIP}.

Next we define the $C([0,1]^3)$-valued process
\begin{align} \label{eq:What}
  \hat{\mathbb{W}}_n(s,t,u)
  = \hat {\mathbb{V}}_n(s,t,u) - s \hat{\mathbb{V}}_n(1,t,u) \, ,
  \qquad s,t,u \in [0,1] \, ,
\end{align}
then the convergence in \eqref{eq:WIP} and the continuous mapping theorem yield
\begin{align} \label{eq:W-convergence}
  \hat {\mathbb{W}}_n \rightsquigarrow \mathbb{W}
\end{align}
in $C([0,1]^3)$, where $\mathbb{W}$ is centred Gaussian defined by
$
  \mathbb{W}(s,t,u)= \mathbb{V}(s,t,u) - s \mathbb{V}(1,t,u)
$
with covariance operator given by \eqref{eq:W-cov}.
Finally, recall the definition of the process $(\hat{\mathbb{U}}_{n}\colon n\in\mathbb{N})$  in \eqref{un}
and note   that, in contrast to $\hat {\mathbb{W}}_n$, this process is not centred.
Consequently, if    $d_\infty = 0$,  we have $\sqrt{n}\, \mathbb{U}_n = \hat {\mathbb{W}}_n$
and  the convergence in \eqref{eq:W-convergence} and the
continuous mapping theorem directly yield \eqref{teil1}.

If $d_\infty > 0$   assertion \eqref{eq:D-limit} is a consequence of the weak convergence in \eqref{eq:W-convergence} and
 Theorem B.1 in the online supplement of \cite{dette2018} and also of the results in \cite{carcamo2019}.

\subsection{Proof of Theorem~ \ref{thrm:bootstrap_test_cpclass} and \ref{thrm:bootstrap_test_cp}} \label{65}

{\bf Proof of Theorem \ref{thrm:bootstrap_test_cpclass}.}
It can be shown that the bootstrap processes in \eqref{eq:bootProcess-cp} can be written
\begin{align*}
\begin{split}
\hat{B}_n^{(r)}(s,t,u) =& \frac{1}{\sqrt{n}} \sum_{k=1}^{\lfloor sn \rfloor}
\frac{1}{\sqrt{l}} \Big( \sum_{j=k}^{k+l-1} \tilde{Y}_{n,j}(t,u)
- \frac{l}{n} \sum_{j=1}^n \tilde{Y}_{n,j}(t,u) \Big) \xi_k^{(r)} \\
&+ \sqrt{n}\Big(s - \frac{\lfloor sn \rfloor}{n} \Big)\frac{1}{\sqrt{l}}
\Big( \sum_{j=\lfloor sn \rfloor +1}^{\lfloor sn \rfloor+l} \tilde{Y}_{n,j}(t,u)
- \frac{l}{n} \sum_{j=1}^n \tilde{Y}_{n,j}(t,u) \Big)
\xi_{\lfloor sn \rfloor +1}^{(r)}  + o_\P(1)
\end{split}
\end{align*}
for $r = 1,\dots, R$ where
\[
\tilde{Y}_{n,j} = \tilde{\eta}_{n,j}^{\check \otimes 2}(t,u) - (\hat{C}_2 - \hat{C}_1)
\mathds{1}\{j > \lfloor \hat{s}n \rfloor \}
\]
for $j=1,\dots,n$ ($n\in\mathbb{N}$) and the array
$(\tilde{\eta}_{n,j}^{\check \otimes 2} \, \colon n\in\N, ~ j = 1,\dots, n)$ satisfies (A1), (A3) and (A4) of
Assumption~2.1 in \cite{dette2018}. The convergence in \eqref{eq:W-convergence} and similar arguments as in the proof of
Theorem~4.3 in the same reference then imply
\begin{align} \label{eq:vec-W}
\big( \hat{\mathbb{W}}_n ,
\hat{\mathbb{W}}_n^{(1)},\dots,\hat{\mathbb{W}}_n^{(R)}\big)
\rightsquigarrow
(\mathbb{W}, \mathbb{W}^{(1)},\dots, \mathbb{W}^{(R)})
\end{align}
in $C([0,1]^3)^{R+1}$ as $n \to \infty$ where the process $\hat{\mathbb{W}}_n$ is
defined by \eqref{eq:What}, the bootstrap counterparts
$\hat{\mathbb{W}}_{n}^{(1)},\dots,\hat{\mathbb{W}}_{n}^{(R)}$ are defined by
\eqref{eq:What-boot} and the random variables
$\mathbb{W}^{(1)},\dots,\mathbb{W}^{(R)}$ are independent copies
of $\mathbb{W}$ which is defined by its covariance operator \eqref{eq:W-cov}.

If $d_\infty = 0$, the continuous mapping theorem directly implies
\begin{align*}
  \big( \hat{\mathbb{M}}_n ,
  \check T_{n}^{(1)},\dots,\check T_{n}^{(R)}\big)
  \stackrel{\mathcal{D}}{\longrightarrow}
  (\check T, \check T^{(1)},\dots, \check T^{(R)})
\end{align*}
in $\R^{R+1}$ as $n \to \infty$ where the statistic $\hat{\mathbb{M}}_n$ is
defined by \eqref{eq:M-hat}, the bootstrap statistics
$\check T_{n}^{(1)},\dots,\check T_{n}^{(R)}$ are defined by
\eqref{eq:checkT} and the random variables
$\check T^{(1)},\dots,\check T^{(R)}$ are independent copies
of the random variable $\check T$ defined by \eqref{teil1}. Now the same arguments as in the discussion starting from equation \eqref{eq:vec-T} imply the assertions made in Theorem~\ref{thrm:bootstrap_test_cpclass}.

\medskip

 {\bf Proof of   Theorem~\ref{thrm:bootstrap_test_cp}.}
  We first mention that it follows by similar arguments as given in the proof of Theorem~4.2 in \cite{dette2018} that the estimator of the unknown change location defined by \eqref{cpEstimator} satisfies
$$|\hat{s}-s^*| = O_{\mathbb{P}}(n^{-1})$$
whenever $d_\infty >0$. Whenever $d_\infty =0$, suppose that
the estimate $\hat{s}$ converges weakly to a $[\vartheta,1-\vartheta]$-valued
random variable which is denoted by $s_{\max}$. Then, if $d_\infty > 0$, the convergence in \eqref{eq:D-limit} and Slutsky's theorem yield
\begin{align} \label{eq:D-limit2}
\sqrt{n}\big( \hat{d}_\infty - d_\infty \big) \stackrel{\cal D}{\longrightarrow}
D (\mathcal{E})= {\tilde D (\mathcal{E}) }/[{s^*(1-s^*)}] \, ,
\end{align}
where $\tilde D (\mathcal{E})$ is the same as in \eqref{eq:D-limit} and the
statistic $\hat{d}_\infty$ is defined by \eqref{eq:statistic-cp}.

The same arguments as in the proof of Theorem~3.6 in \cite{dette2018} again
yield that the estimators of the extremal sets defined by \eqref{eq:estimatedSets-cp}
are consistent. The convergence in \eqref{eq:vec-W} and similar
arguments as in the proof of Theorem~4.4 in the same reference then yield
\begin{align} \label{eq:vec-convergence-cp}
\big( \sqrt{n} ~ (\hat{d}_\infty - d_\infty) ,~
\check K_{n}^{(1)},\dots,\check K_{n}^{(R)}\big)
\stackrel{\mathcal{D}}{\longrightarrow}
(D(\mathcal{E}),~ D^{(1)}(\mathcal{E}),\dots, D^{(R)}(\mathcal{E}))
\end{align}
in $\R^{R+1}$ as $n \to \infty$ where the bootstrap statistics
$\check K_{n}^{(1)},\dots,\check K_{n}^{(R)}$ are defined by
\eqref{eq:bootStat-cp} and the random variables
$D^{(1)}(\mathcal{E}),\dots,D^{(R)}(\mathcal{E})$ are independent copies
of $D(\mathcal{E})$ which is defined by \eqref{eq:D-limit2}. The
convergence in the preceding equation holds true under the null and the alternative hypothesis and now the same arguments as in the discussion starting from equation \eqref{eq:vec-convergence} imply the assertions made in Theorem~\ref{thrm:bootstrap_test_cp}.

\bigskip
\medskip

\noindent
{\bf Acknowledgements}
This research was partially supported by the Collaborative Research Center `Statistical modeling
of nonlinear dynamic processes' ({\it Sonderforschungsbereich 823, Teilprojekt A1, C1})
and the Research Training Group `High-dimensional phenomena in probability - fluctuations and
discontinuity' ({\it RTG 2131}).
The authors are grateful to Christina Stoehr for sending us the results of
 \cite{stoehr2019} and to  Martina
Stein, who typed  parts  of this manuscript with considerable technical expertise.

\end{document}